\documentclass[letterpaper]{article}
\usepackage{authblk}
\usepackage{url}
\usepackage[margin=1in,footskip=0.25in]{geometry}
\usepackage[labelfont=bf]{caption}
\usepackage{graphicx}
\usepackage{amsmath,amssymb,amsfonts,amsmath}
\usepackage{graphicx}
\usepackage{subfigure}
\usepackage{caption}
\usepackage{tikz}
\usepackage{float}
\usepackage[dvipdf]{epsfig}
\usepackage{float}
\usepackage{lipsum}
\usepackage{wrapfig}
\usepackage{cite}
\usepackage{color}
\usepackage{easylist}
\usepackage{lineno}
\usepackage{booktabs}
\usepackage{hyperref}
\usepackage{cleveref}
\usepackage{multirow}
\usepackage{hhline}
\usepackage{soul}
\usepackage{blindtext}
\usepackage{threeparttable}
\usepackage[bbgreekl]{mathbbol}

\def \L2{L^2}

\newcommand{\eqn}[1]{
	\begin{eqnarray}
	#1
	\end{eqnarray}
}

\newcommand{\tb}[1]{
	{\bf #1}
}

\newcommand{\norm}[1]{
	\lVert #1\rVert
}

\newcommand{\ip}[1]{
	\langle #1\rangle
}

\title{The Smooth Forcing Extension Method: A High-Order Technique for Solving Elliptic Equations on Complex Domains}

\author[1]{Saad Qadeer} 
\author[2]{Boyce E.~Griffith} 
\affil[1]{Department of Mathematics, University of North Carolina, Chapel Hill, NC, USA}
\affil[2]{Departments of Mathematics, Applied Physical Sciences, and Biomedical Engineering, University of North Carolina, Chapel Hill, NC, USA}
\affil[ ]{{\tt saadq@email.unc.edu} and {\tt boyceg@email.unc.edu}}

\date{}

\begin{document}

\maketitle

\begin{abstract}
High-order numerical methods for solving elliptic equations over arbitrary domains typically require specialized machinery, such as high-quality conforming grids for finite elements method, and quadrature rules for boundary integral methods. These tools make it difficult to apply these techniques to higher dimensions. In contrast, fixed Cartesian grid methods, such as the immersed boundary (IB) method, are easy to apply and generalize, but typically are low-order accurate. In this study, we introduce the Smooth Forcing Extension (SFE) method, a fixed Cartesian grid technique that builds on the insights of the IB method, and allows one to obtain arbitrary orders of accuracy. Our approach relies on a novel Fourier continuation method to compute extensions of the inhomogeneous terms to any desired regularity. This is combined with the highly accurate Non-Uniform Fast Fourier Transform for interpolation operations to yield a fast and robust method. Numerical tests confirm that the technique performs precisely as expected on one-dimensional test problems. In higher dimensions, the performance is even better, in some cases yielding sub-geometric convergence. We also demonstrate how this technique can be applied to solving parabolic problems and for computing the eigenvalues of elliptic operators on general domains, in the process illustrating its stability and amenability to generalization.         

\end{abstract}



\section{Introduction}\label{secintro}
A long-standing challenge in the numerical study of elliptic partial differential equations is the development of high-order methods for arbitrary domains. Over the years, various approaches have been proposed and extensively analyzed, refined, and applied to problems from diverse settings. The finite element and boundary integral methods provide elegant formulations of the problem and yield powerful solvers. However, they require additional tools and machinery that limit their scope and hinder their generalization to higher dimensions. The finite element method, for instance, requires a high quality mesh \cite{shewchuk2002good}; in the case of moving boundaries, generating such conforming grids at each time-step can be computationally infeasible. Similarly, boundary integral methods require specialized quadrature rules to resolve the nearly singular kernels \cite{helsing2008evaluation,beale2001method,fryklund2019integral}; while they perform impressively in two dimensions, it is unclear how these tools optimally extend to three dimensions.

At the other end of the spectrum are fixed Cartesian grid methods. Broadly speaking, these techniques embed the physical domain in a simpler computational domain and solve the appropriately modified equations on a non-conformal structured mesh. The immersed boundary (IB) method was introduced by Peskin \cite{peskin1972flow,peskin1977numerical} for simulating fluid flow around immersed elastic bodies. Applied to elliptic problems, it operates by extending the inhomogeneous terms trivially to the computational domain. Any resulting discontinuities in the derivative of the solution are accounted for by the addition of singular terms to the equations. These compactly supported spreading terms act as Lagrange multipliers to enforce the boundary conditions. This formulation is combined with a finite difference discretization, with regularized delta functions used to numerically handle the spreading and interpolation operations. The technique yields first-order accuracy but possesses the key advantages of straightforward generalizability to higher dimensions, the ability to handle moving boundaries, and amenability to adaptive mesh refinement \cite{griffith2007adaptive}. As a result, this method has been successfully used for a wide range of problems \cite{bhalla2013unified,griffith2012immersed,huang2012three,kou2015fully,seol2016immersed}.


An alternative approach to achieving higher-order accuracy in the IB method is the Immersed Boundary Smooth Extension (IBSE) method \cite{stein2016immersed,stein2017immersed}. The key insight in the design of this technique is that the bottleneck in achieving high-order accuracy is the smoothness of the posited solution on the computational domain. The method poses a high-order PDE outside the physical domain, with boundary conditions matching the exact solution. This problem is solved using the IB method, and the solution extension in turn is used to supply the necessary extension to the forcing. The modifications indeed lead to high-order accuracy; however, its use of regularized delta functions limits its efficacy. More damagingly, its introduction of a high-order equation yields an unwieldy structure that, in some cases, may lead to ill-conditioning and instabilities.

In this paper, we present a new fixed Cartesian grid method for solving elliptic problems by further developing the insights introduced in the IBSE method. Our technique eschews directly solving for a smooth extension to the unknown solution and instead uses a novel Fourier extension method to extend the forcing. This approach demonstrably resolves the ``mountain-in-fog'' problem and can be used to compute the extension to any desired regularity \cite{boyd2002comparison}. Coupled with highly accurate inversion and interpolation procedures using the Fast Fourier Transform (FFT), this leads to a rapid, robust, and highly accurate technique for solving elliptic equations. A hallmark of our approach is its simplicity, which allows it to be used for complex domains (including those with sharp corners) in any number of dimensions. Moreover, our method possesses strong stability properties that, as we shall demonstrate, allow it to be extended to parabolic problems and be used to compute the eigenvalues of elliptic operators on arbitrary domains.

Our technique is also amenable to other discretization approaches. For instance, one could employ a finite difference discretization based on a uniform grid and make use of fast iterative solvers, notably multigrid methods, instead of the Fourier solver. It must be ensured, however, that the corresponding interpolation procedures are not based on high-order polynomials to avoid the instabilities associated with uniform grids.

Another alternative is to use a tensor product of one-dimensional Chebyshev--Lobatto grids. This technique has the advantage of allowing rapid FFT-based transforms to representations in terms of Chebyshev polynomials \cite{trefethen2013approximation}. These polynomials are inexpensive to differentiate and lend themselves to fast and accurate evaluation at off-grid points by formulas based on the Clenshaw recurrence formula \cite{clenshaw1955note}. Thus, this choice enables interpolation to be performed efficiently. However, this comes at the cost of being unable to use the fast finite difference solvers. Instead, we are required to use a Galerkin formulation that leads to dense stiffness matrices, making this approach somewhat prohibitive.

Prior efforts devoted to improving upon the IB method have also led to the development of widely used techniques such as the Immersed Interface method \cite{leveque1994immersed,li2006immersed}, the Ghost Fluid method \cite{fedkiw1999non,marques2011correction}, and the Active Penalty method \cite{shirokoff2015sharp}.  Another recent technique that uses an approach similar to that described herein is the Fourier Continuation Alternate Direction (FC-AD) Implicit method pioneered by Bruno and Lyon \cite{bruno2010high,lyon2010high}. This approach relies on the ADI procedure to reduce an evolution equation to a sequence of one-dimensional elliptic problems, which are extended by a highly accurate Fourier Continuation routine to the appropriate computational domains. Also of note is the recent work on smooth selection embedding, which attempts to solve the extension problem by formulating it as a Sobolev norm optimization problem \cite{agress2018novel,agress2019smooth}. In general, Fourier continuation methods have a rich history; see \cite{boyd2002comparison} for an exhaustive review. Yet another exciting contribution is the partition of unity extension approach developed in the context of boundary integral methods and applied to heat and fluid flow problems \cite{fryklund2019integral,fryklund2018partition,klinteberg2019fast}.


\section{The Smooth Forcing Extension Method}\label{secsfe}  

\subsection{Mathematical Formulation}\label{secmath}
We describe our method by outlining its use for the problem 
\eqn{
	\left\{ \begin{matrix}
		\mathcal{L} u = f, &  \text{on } \Omega,  \\
		u = g, &   \text{on } \partial \Omega.  \\
	\end{matrix} \right.
	\label{model}
}
Here, $\Omega$ is an arbitrary bounded domain in $\mathbb{R}^d$ or $\mathbb{T}^d$, $\mathcal{L}$ is an elliptic operator, and $f \in C^{\infty}(\Omega)$; we shall frequently refer to $f$ as the forcing. For clarity of exposition, we have restricted ourselves to a Dirichlet problem for now; we shall later show that our technique can easily handle all types of boundary conditions.

We begin by embedding $\Omega$ in a computational domain $C$ and defining the extension region $E = C - \overline{\Omega}$. In addition, for $k \geq 0$, let $T_k^*$ denote the evaluation operator for the first $k$ normal derivatives at the boundary. Note then that $S^* := T_0^*$ is simply the interpolation operator.

An important step in the development of the IBSE method was the observation that, in principle, $u$ can be extended smoothly to $C$. Contrary to how this technique proceeds, however, we shall not explicitly solve for the extension to the solution. Instead, we further note that any additional forcing induced by the extended solution must remain restricted to $E$, as $\mathcal{L}$ is a local operator. Thus, we can search for the extension to the forcing in a space of functions supported on $E$.

Let $\{\phi_j\}_{1 \leq j \leq J}$ be a family of smooth functions on $C$; this shall serve as the basis of the space in which we shall look for the extension to the forcing. Set $h = \sum_{j = 1}^J c_j\phi_j$ and consider the extended problem   
\eqn{
	\left\{ \begin{matrix}
		\mathcal{L} u_\text{e} = \chi_{\Omega}f + \chi_E h, &  \text{on } C,  \\
		S^*u_\text{e} = g, &  \text{on } \partial \Omega.   \\
	\end{matrix} \right.
	\label{sfe}
}
Here, $\chi_A$ denotes the characteristic function for a set $A$. Assuming that $\mathcal{L}$ is invertible on $C$, we obtain
\eqn{
	u_\text{e} = \mathcal{L}^{-1}(\chi_{\Omega}f) + \mathcal{L}^{-1}(\chi_E h) = \mathcal{L}^{-1}(\chi_{\Omega}f) + \sum_{j = 1}^J c_j \mathcal{L}^{-1}(\chi_E \phi_j). \label{Ueeq}
}
Applying $S^*$ throughout and using $S^*u_\text{e} = g$ yields
\eqn{
	\sum_{j = 1}^J c_j S^*\mathcal{L}^{-1}(\chi_E \phi_j) = g - S^*\mathcal{L}^{-1}(\chi_{\Omega}f). \label{Seq}
}
Next, observe that if $u_\text{e}$ is sufficiently smooth on $C$, some of this regularity would be inherited by the extended forcing $f_\text{e} = \chi_{\Omega}f + \chi_E h$. This condition can be enforced by requiring that 
\eqn{
	\sum_{j = 1}^J c_j T_k^*\phi_j = T_k^*h = T_k^*f, \label{Tkeq}
}
for some $k \geq 0$. It follows from elliptic regularity theory that if $f_\text{e} \in C^k(C)$ and $\mathcal{L}$ is of order $l$, then $u_\text{e} \in C^{k+l}(C)$ \cite{evans2010partial}. Note, however, that $f$ may not be known outside $\Omega$, or, its analytic continuation may contain singularities, so $T_k^*f$ may be ill-defined. To remedy this, we may instead enforce $T_k^*h = T_k^*f_e$, which reduces to
\eqn{
	\sum_{j = 1}^J c_j T_k^*(\chi_{\Omega}\phi_j) = T_k^*(\chi_{\Omega}f). \label{Tkeq2}
}
Taken together, equations (\ref{Seq}) and (\ref{Tkeq}) (or (\ref{Tkeq2})) prescribe the conditions that must be met to yield a problem with a sufficiently smooth solution on $C$. These conditions are to be satisfied at the boundary so a discretization $\tb{s} = (s_i)_{1 \leq i \leq n_\text{b}}$ of $\partial \Omega$ would result in a linear system of size $n_{b}(k+2) \times J$. Thus, by choosing $J$ so that the system is square (or under-determined), we can solve for the coefficients $\{c_j\}$ (in the minimum norm sense) and obtain the extended solution $u_\text{e}$ using (\ref{Ueeq}).

\subsection{Implementation Details}\label{secimpl}
After outlining the basic ideas behind our method, we shall now discuss some details regarding its implementation that allow us to harness its full accuracy, efficiency and robustness.

Our approach to discretization is aimed at making full use of Fourier-based techniques. The computational domain $C$ is taken as the $d$-dimensional periodic box $\mathbb{T}^d$, with equal-sized grid cells, and the extension functions are chosen as the trigonometric polynomials $\{e^{i{\bf j}\cdot {\bf x}}\}$. A significant advantage of this approach is the simplicity, speed and accuracy of inverting the differential operator $\mathcal{L}$: the derivative $\partial_{x_1}^{\alpha_1}\partial_{x_2}^{\alpha_2}\hdots \partial_{x_d}^{\alpha_d}$ is replaced by its symbol $(i\xi_1)^{\alpha_1}(i\xi_2)^{\alpha_2}\hdots (i\xi_d)^{\alpha_d}$ in Fourier space, so that implementing $\mathcal{L}^{-1}$ reduces to a pair of FFTs and a term-wise algebraic solve. In addition, we can take advantage of Non-Uniform FFT (NUFFT) algorithms to discretize $T_k^*$ extremely accurately and efficiently \cite{greengard2004accelerating,lee2005type}. More precisely, given a function $f$ on $C$, a discretization $\tb{s} = (s_i)_{1 \leq i \leq n_\text{b}}$ of the boundary $\partial \Omega$ and unit normal vectors $\{\tb{n}_i\}_{1\leq i \leq n_\text{b}}$ at the respective nodes, we compute 
\eqn{
	T_k^*f = \left(\begin{matrix}
		f(\tb{s}) & D_{\tb{n}}f(\tb{s}) & D^2_{\tb{n}}f(\tb{s}) & \hdots & D^k_{\tb{n}}f(\tb{s}) \\ 
	\end{matrix}\right)^T	\label{Tks}
}
where $D^l_{\tb{n}}f(\tb{s}) = \left( D^l_{\tb{n}_i}f(s_i)\right)_{1 \leq i \leq n_\text{b}}$ consists of the $l$th directional derivatives of $f$ at all the boundary nodes, in the direction of the corresponding normal vectors. If $\tb{n}_i = (a^{(i)}_1 \quad a^{(i)}_2 \quad \hdots \quad a^{(i)}_d )^T$, these normal derivatives are given by 
\eqn{
	D^l_{\tb{n}_i}f(s_i) = \left(\sum_{j = 1}^d a^{(i)}_j \partial_{x_j}\right)^l f(s_i) = \sum_{\alpha_1 + \hdots + \alpha_d = l} {l \choose \alpha_1,\hdots,\alpha_d} \prod_{m = 1}^d \left(a^{(i)}_m \partial_{x_m} \right)^{\alpha_m} f(s_i). \label{dirder}	
}  
The partial derivatives of $f$ can be computed in Fourier space and evaluated at the boundary nodes by using NUFFT.      

An issue with this approach is that we may lose the invertibility of $\mathcal{L}$ that we made use of earlier since $\mathbb{T}^d$ has no boundaries (so we cannot impose additional boundary conditions on $\partial C$). This would be the case, in particular, if $\mathcal{L} = \Delta$. We specify the recipe for this example because of its ubiquity, although a similar procedure can be followed for any self-adjoint $\mathcal{L}$. Decompose 
\eqn{
	u_\text{e} = U + u_0, \label{Udefn}
}
where $U = |C|^{-1} \int_C u_\text{e} \ dx$, so that $\int_C u_0 \ dx = 0$ and $\Delta u_\text{e} = \Delta u_0$. One can then replace $\mathcal{L}^{-1}$ in (\ref{Seq}) by the ``zero-mean'' inverse $\mathcal{A}$ of the Laplacian to obtain
\eqn{
	\sum_{j = 1}^J c_j S^*\mathcal{A}(\chi_E \phi_j) + U = g - S^*\mathcal{A}(\chi_{\Omega}f). \label{Seq2}
}

In addition, we average the first equation in (\ref{sfe}) over $C$ and use $\int_C \Delta u_\text{e} \ dx = \int_C \Delta u_0 \ dx = 0$ to obtain the additional equation
\eqn{
	\sum_{j = 1}^J c_j \int_C  \chi_E \phi_j \ dx = -\int_C \chi_{\Omega} f \ dx. \label{Avgeq}
}

As above, equations (\ref{Seq2}) and (\ref{Avgeq}) can be complemented with the regularity constraints (\ref{Tkeq}) (or (\ref{Tkeq2})) to form a system of size $(n_{b}(k+2)+1) \times (J+1)$ and can be used to solve for $u_\text{e}$. In the case the system is under-determined, the solutions are not unique. We then use the Moore--Penrose pseudo-inverse to find the minimum norm solution \cite{demmel1997applied}.

\section{Numerical Results} \label{numres}
In order to demonstrate the effectiveness of our algorithm, we shall present results from a variety of contexts. We begin by solving some simple problems in one dimension to further elucidate its implementation and assess its performance, before moving to test problems in two dimensions. 

\subsection{The Extension Algorithm}
As a preliminary test, we investigate the extension routine: given a function $f$ on $\Omega$, we extend it to $f_\text{e}$ on the computational domain $C = \mathbb{T}^d$ such that $f_\text{e}|_{\Omega} = f$. A highly desirable property of an extension algorithm is that it circumvent the ``mountain-in-fog'' problem \cite{boyd2002comparison}. This refers to the pitfall that an algorithm chooses the analytic continuation of $f$ that may contain singularities in the extension region.

For a simple one-dimensional example of such a function, let $\Omega = (2,5)$ and $f(x) = 1/(x-1)$. We define the extension by 
\eqn{
	f_\text{e}(x) = \chi_{\Omega}f(x) + \chi_E\sum_{j = -J}^J c_je^{ijx}, \label{fext}
}
where $E = \mathbb{T}-\bar{\Omega}$; the reality conditions $c_{-j} = c_j^*$ for all $j$ imply that we effectively have $(2J+1)$ real degrees of freedom. The only conditions on the $\{c_j\}$ are the regularity constraints (\ref{Tkeq}) that provide matching conditions for the values and first $k$ derivatives at $\partial \Omega$. For a given $k$, we choose $J = k+1$ to obtain a system of size $(2k+2) \times (2k+3)$, which is solved (in the minimum norm sense) to yield the $k$-regular extension $f_\text{e}^{[k]}$. In this computation, the right hand side of (\ref{Tkeq}), $T_k^*f$, is calculated exactly since $f$ is known in closed-form. If $f$ was only known at the grid points in $\Omega$, we would have instead used regularity conditions of the form given in equation (\ref{Tkeq2}).

\def \sclbs {0.35}
\begin{figure}[tbph]
	\centering
	\subfigure[]
	{\includegraphics[scale=\sclbs]{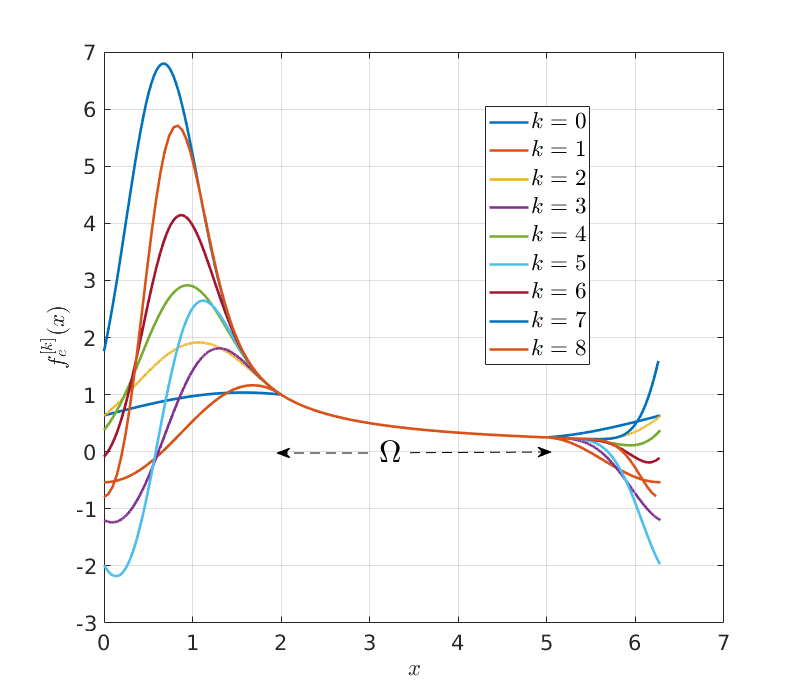}
		\label{fextcurv1}
	}
	\subfigure[]
	{\includegraphics[scale=\sclbs]{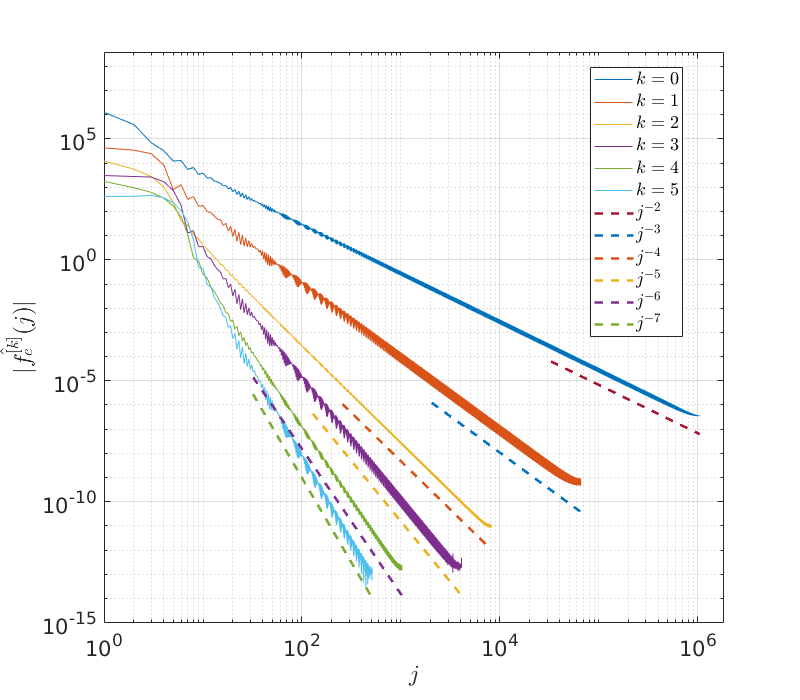}
		\label{fectdec1}
	}
	\caption{(a) The extensions $f_\text{e}^{[k]}$ for different $k$, for $f(x) = 1/(x-1)$ and $\Omega = (2,5)$. Note that the functions approach the singularity at $x = 1$ for higher $k$, since that is a feature of the analytic continuation, but still remain well-behaved. (b) The Fourier coefficients of the extensions shown in (a) decay like $O(j^{-(k+2)})$. Results for higher $k$ are omitted for clarity. }\label{fext1}
\end{figure}

If $f^{[k]}_e \in C^k(\mathbb{T})$, its Fourier coefficients should decay asymptotically as $O(j^{-(k+2)})$. The extensions for different values of $k$ are shown in Figure \ref{fext1} along with the decay of their Fourier coefficients. The decay rates are as expected, indicating the effectiveness of our extension technique for arbitrary $k$. Higher values of $k$ also exhibit the same trend but are omitted from the plot. A notable point is that, while the extensions in Figure \ref{fextcurv1} appear to approach the singularity at $x = 1$ for increasing values of $k$, they avoid the ``mountain-in-fog'' problem. This is primarily because we only use the boundary data and search for the extension in a low-dimensional space. Traditional techniques extrapolate the function after sampling it in the interior of $\Omega$ and, as a result, are more likely to mimic the pathological behavior. 

\subsection{Poisson Equation in One Dimension} 
Next, we consider the Poisson equation
\eqn{
	\left\{ \begin{matrix}
		u_{xx} = 1/(x-1), &  \text{on } \Omega = (2,5), \\
		u(2) = u_2, &    u(5) = u_5.  \\
	\end{matrix} \right.
	\label{lap1Dprob}
}
The exact solution to this problem can be used for comparison against the numerical solutions. To calculate the forcing extensions, we also impose the boundary conditions (\ref{Seq2}) and averaging condition (\ref{Avgeq}); as a result, the particular extensions in this case will be different from those shown in Figure \ref{fextcurv1} while possessing the same regularity. In fact, we can also compute a forcing extension without imposing any smoothness requirements; we refer to this as $k = -1$; the resulting extension has a jump discontinuity at $\partial \Omega$. Finally, for a given $k$, we set $J = k+2$ to yield an under-determined system of size $(2k+5) \times (2k+6)$.

\def \sclas {0.5}
\begin{figure}[tbph]
	\centering
	{\includegraphics[scale = \sclas]{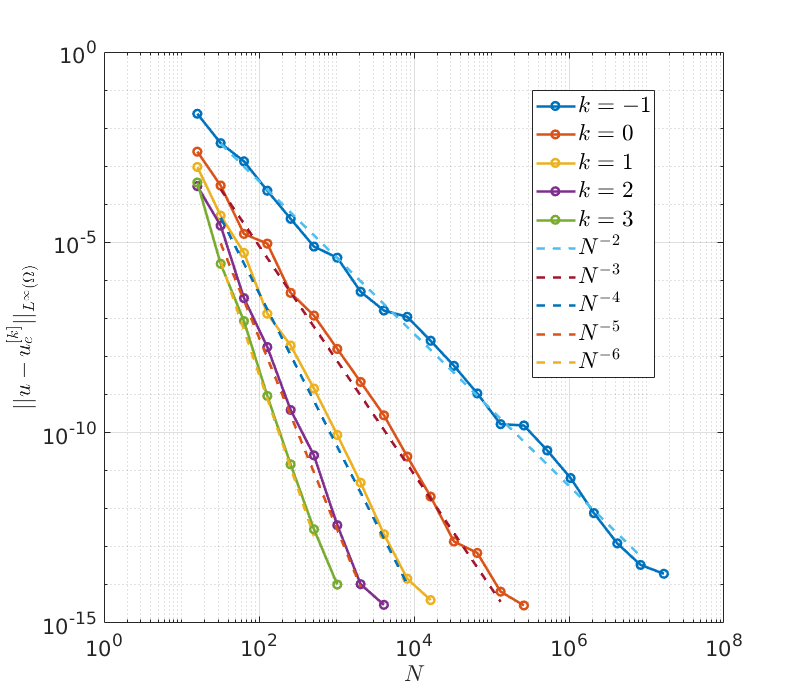}}
	\caption{Errors in the $L^{\infty}$ norm in the numerically computed solutions for problem (\ref{lap1Dprob}). The solution derived from forcing with $k$ continuous derivatives on $C$ can be seen to converge at rate $O(N^{-(k+3)})$, for arbitrary $k$.} \label{lap1dinf}
	
\end{figure}

Since the operator $\mathcal{L} = \partial_x^2$ is second-order, by the earlier discussion, a $k$-regular extension should yield a solution $u_\text{e}^{[k]} \in C^{k+2}(C)$. As a result, we expect convergence in the $L^2$ norm at rate $O(N^{-(k+3)})$, where $N$ is the number of grid points. Figure \ref{lap1dinf} shows the results with $u_2 = 1$ and $u_5 = -1$ in the $L^{\infty}$ norm. We employ this norm as it bounds the $L^2$ norm while also allowing us to assess convergence at points close to the boundary. It can be seen that the errors indeed converge at the desired rates and achieve 13 digits of accuracy in all the cases.

\begin{figure}[tbph]
	\centering
	{\includegraphics[scale = \sclas]{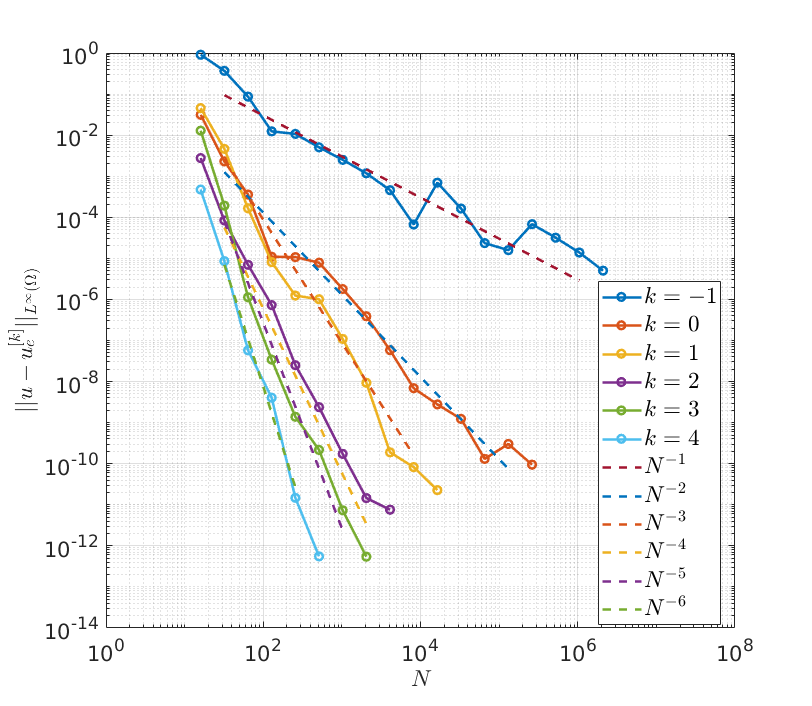}}
	\caption{$L^{\infty}$ errors in the numerically computed solutions for the 1D Poisson equation with mixed boundary conditions. The solutions converge to the true solution at $O(N^{-(k+2)})$; the reduction in order is due to the loss of accuracy while imposing the Neumann condition.} \label{lap1dinfmix}
	
\end{figure}

Small modifications in our method allow us to handle different boundary conditions. Consider the same problem as (\ref{lap1Dprob}) with $u_x(2) = u_2$ and $u(5) = u_5$. The Neumann condition is imposed by changing (\ref{Seq2}) to 
\eqn{
	\sum_{j = 1}^J c_j \left[T_1^*\mathcal{A}(\chi_E \phi_j)\right]_{x = 2} &=& u_2 - \left[T_1^*\mathcal{A}(\chi_{\Omega}f)\right]_{x = 2}. \label{Seq2N}
}
The Dirichlet condition at $x = 5$ is imposed in exactly the same manner. Note that the mean correction $U$ does not appear in (\ref{Seq2N}) since it vanishes upon differentiation. The results, for $u_2 = 1$ and $u_5 = -1$ with the same number of extension functions as above, are shown in Figure \ref{lap1dinfmix}. Observe that the errors decay at $O(N^{-(k+2)})$. The reduction in order is due to the fact that the accuracy of the derivative--interpolation operator $T_j^*$ decreases with increasing $j$. Thus, applying condition (\ref{Seq2N}) introduces a bottleneck, which is reflected in the error decays of the solutions.

\subsection{Heat Equation in One Dimension}
As our final one-dimensional example, we show how to adapt this method to solve the heat equation. Consider the more general time-dependent problem 
\eqn{
	\left\{ \begin{matrix}
		u_t - \mathcal{L} u = f(t,x), &  \text{for }x \in \Omega,  \\
		u(0,x) = u_0(x), &  \text{for }x \in \Omega,  \\
		u(t,x) = g(t,x),  & \text{for }x \in  \partial \Omega, \ t > 0.   \\
	\end{matrix} \right.
	\label{tdepeqn}	
}

\def \sclam{0.2}
\begin{figure}[tbph]
	\centering
	{\includegraphics[scale = \sclam]{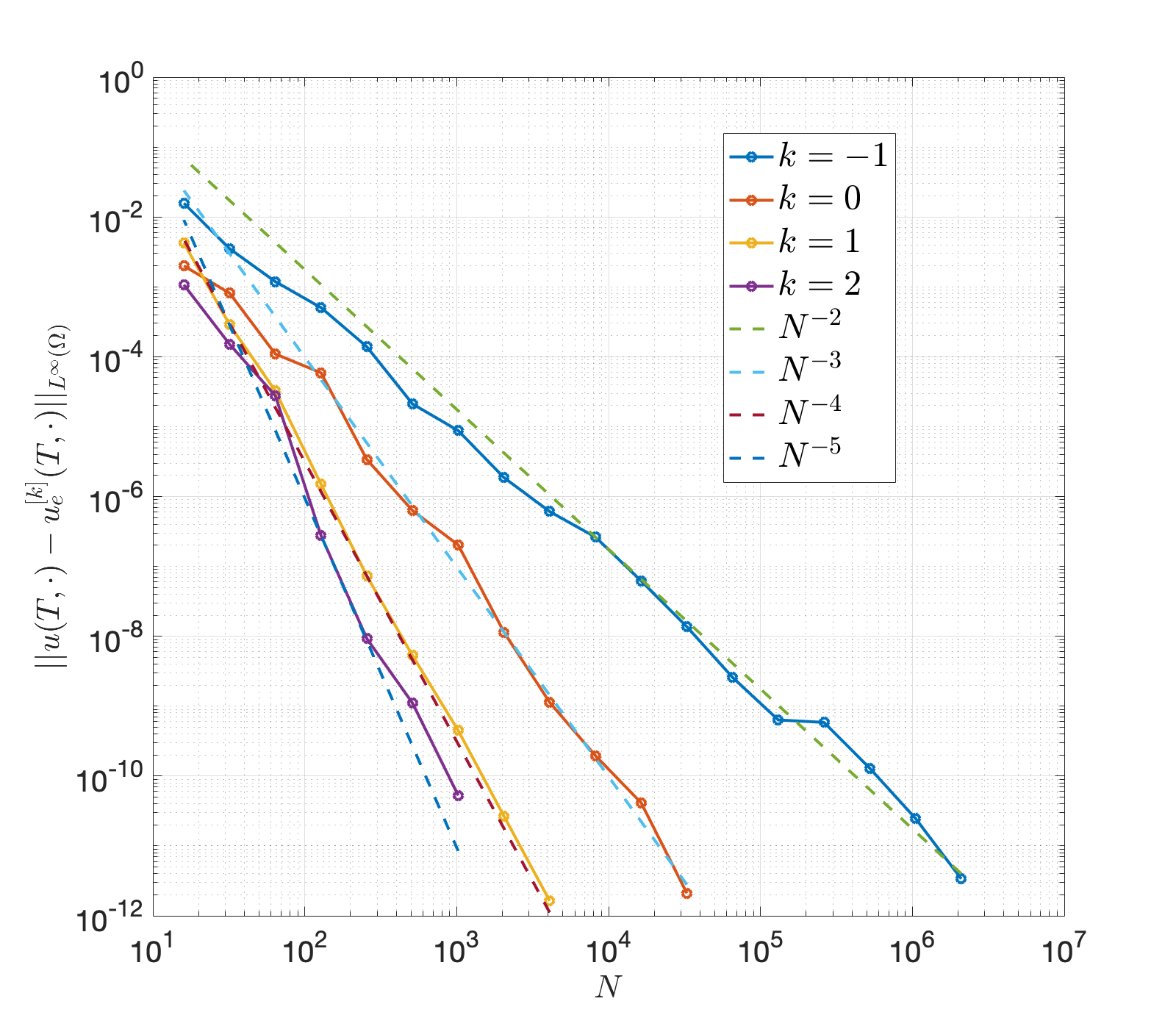}}
	\caption{$L^{\infty}$ errors for problem (\ref{heat1dprob}). The solutions can be seen to converge more or less as $O(N^{-(k+3)})$.} \label{heat1dinf}
	
\end{figure}

We employ the following iteration scheme, obtained from the four-step Backward Differentiation Formula (BDF-4), to discretize the time derivative:
\eqn{
	\left(\mathbb{I} - \frac{12\Delta t}{25} \mathcal{L}\right)u^{n+1} = \frac{12\Delta t f^{n+1} + 48u^{n} - 36u^{n-1} + 16u^{n-2} - 3u^{n-3}}{25},	\label{bdf4}
}
Setting
\eqn{
	\hat{ \mathcal{L}} = \mathbb{I} - \frac{12\Delta t}{25} \mathcal{L}, \qquad F^{n+1} = \frac{12\Delta t f^{n+1} + 48u^{n} - 36u^{n-1} + 16u^{n-2} - 3u^{n-3}}{25}	\label{Lhat}
}
allows us to write (\ref{bdf4}) as
\eqn{
	\hat{\mathcal{{L}}}u^{n+1} = F^{n+1} \label{bdf41}	
}
with corresponding boundary conditions $u^{n+1}(x) = g((n+1)\Delta t,x)$, for $x \in \partial \Omega$. This formulation lends itself naturally to the extension technique described earlier. The computationally intensive task of building the matrices corresponding to equations (\ref{Seq}) and (\ref{Tkeq}) needs to be performed just once (for a specified $\Delta t$) for the entirety of a simulation. 

Observe that $F^{n+1}$ requires solutions over four previous time-steps. Hence, for the first three iterations of this technique, we need either additional initial conditions (going back three time-steps) or we need to couple the method with a single-step time integrator for jump-starting the algorithm. In the results presented here, we have followed the second approach since it requires fewer inputs and is more broadly applicable. We opt for the Backward Euler method due to its simplicity and A-stability. However, its use introduces a time discretization error in addition to the $O(\Delta t^4)$ error from the BDF-4 scheme. To ensure that the errors decay at a comparable rate, we further divide each of the first three time intervals into $M$ parts. More precisely, we discretize (\ref{tdepeqn}) as
\eqn{
	\left(\mathbb{I} - \left(\frac{\Delta t}{M}\right)\mathcal{L}\right)u^{(l+1)/M} = \left(\frac{\Delta t}{M}\right) f^{(l+1)/M} + u^{l/M},	\label{bdf1}
}
for $0 \leq l \leq 3M-1$. This procedure yields the approximations $u^n$ for $1 \leq n \leq 3$ that can then be fed into the higher order scheme (\ref{bdf41}). The total time discretization error adds up to $O\left(M\left(\frac{\Delta t}{M}\right)^2 + \Delta t^4\right) = O\left(\frac{\Delta t^2}{M} + \Delta t^4\right)$. As a result, scaling $M$ as $(\Delta t)^{-2}$ ensures that the overall time discretization error is fourth-order.

For an iterative procedure of this form, we find that imposing regularity constraints of the form (\ref{Tkeq}) leads to an unstable system. This happens because building the right-hand sides of (\ref{bdf4}) and (\ref{bdf1}) requires past solutions that possess unphysical extensions outside $\Omega$. Smoothness conditions (\ref{Tkeq}) re-use these values and lead to spurious behavior. Using (\ref{Tkeq2}), meanwhile, avoids this issue altogether as it only makes use of the forcing on $\Omega$, and has the added benefit of not requiring the values of $f(t,x)$ for $x \notin \Omega$.

In the case of the heat equation, we have $\mathcal{L} = \Delta$, so (\ref{bdf41}) reduces to a Helmholtz equation. Figure \ref{heat1dinf} shows the results for $\Omega = (2,5)$ with 
\eqn{
	f(t,x) = \sin(x), \ u_0(x) = e^{\sin(x)}, \  u(t,2) = 1, \ u(t,5) = 0. \label{heat1dprob}
}
The problem was solved up to $T = 1$, with time-step $\Delta t = 2.5 \times 10^{-3}$ and $M = 1$ for all values of $k$ and $N$. As the errors in Figure \ref{heat1dinf} do not appear to plateau, we deduce that time integration errors are negligible, thus allowing for a comprehensive test of the accuracy and stability properties of our algorithm. The asymptotic error decay rates can be seen to be $O(N^{-(k+3)})$.  

\subsection{Poisson Equation in Two Dimensions}
\begin{figure}[tbph]
	\centering
	\subfigure[]
	{\includegraphics[scale=\sclbs]{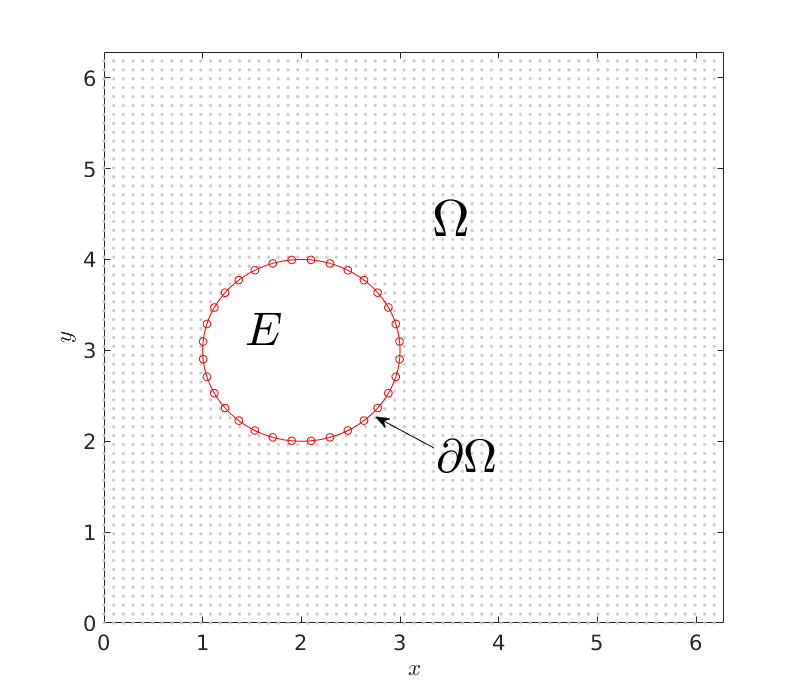}
		\label{discshape}
	}
	\subfigure[]
	{\includegraphics[scale=\sclbs]{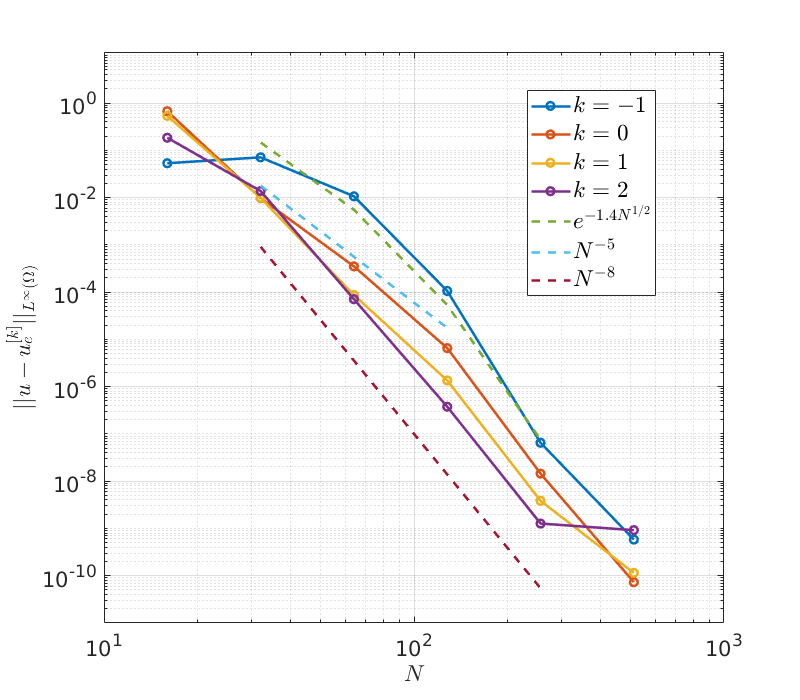}
		\label{lap2dinfdec0}
	}
	\caption{(a) The physical domain $\Omega$ along with its boundary $\partial \Omega$ and extension region $E$. The fixed grid is shown in $\Omega$, along with the boundary nodes, for $N = 2^6$ and $n_\text{b} = 32$. (b) The $L^{\infty}$ errors for a Poisson equation solved on the domain shown in (a). The solutions converge faster than $O(N^{-(k+3)})$ in all cases and sub-geometrically for $k = -1$; a reasonable fit appears to be $O(\exp(-N^{1/2}))$}\label{lap2d0}
\end{figure}

After solving the one-dimensional test problems, we turn our attention towards problems in two dimensions. Let $B_1(2,3)$ be the unit disc centered at $(2,3)$ and define $\Omega = [0,2\pi)^2 - \bar{E}$ (see Figure \ref{discshape}). Consider the Poisson problem  
\eqn{
	\left\{\begin{matrix}
		-\Delta u = 5\sin(x)\cos(y), & \text{on } \Omega, \\
		u = 0, & \text{on } \partial \Omega. 	
	\end{matrix}\right.	\label{lap2dprob1}
}
To solve this, we embed $\Omega$ in $\mathbb{T}^2$, place a uniform grid with $N$ points along each axis, and use the extension family $\{e^{i(j_1x + j_2y)}\}_{-J \leq j_1,j_2 \leq J}$. As in the solution to problem (\ref{lap1Dprob}), we impose the boundary conditions (\ref{Seq2}) and averaging condition (\ref{Avgeq}). The boundary $\partial \Omega$ is discretized by placing a total of $n_\text{b}$ equidistant points on it. We choose $n_\text{b} = \lceil 0.5 N \rceil$ to ensure that the spacing $\Delta s$ between successive boundary nodes is roughly twice the grid spacing $\Delta x$. This ratio has been empirically observed to yield an optimal balance between conditioning and accuracy for fixed grid methods \cite{stein2016immersed,kallemov2016immersed}. In our experiments, we have also found that it leads to superior performance over other choices. 

For a $k$-regular extension, we impose a total of $n_\text{b}(k+2)+1$ constraints. Since we have $(2J+1)^2+1$ degrees of freedom, we set 
\eqn{
	J = \left\lceil \frac{\sqrt{n_\text{b}(k+2)}}{2}-1 \right\rceil \label{Jchoice}
}
to obtain the customary under-determined system. Instead of computing the exact solution to (\ref{lap2dprob1}) by another technique for comparison, we use solutions on successively refined grids to compute the errors. The resulting refinement study, displayed in Figure \ref{lap2dinfdec0}, shows that the technique performs better than expected. The convergence exceeds $O(N^{-(k+3)})$ in all cases and, in particular, is faster than any power of $1/N$ for $k = -1$. The rate is still slower than spectral, and is therefore termed sub-geometric. Indeed, as shown in the plot, $O(\exp(-N^{1/2}))$ models the decay reasonably well.

\def \sclbm{0.143}
\def \sclim{0.038}

\begin{figure}[tbph]
	\centering
	\subfigure[]
	{\begin{tikzpicture}
		[node distance = 10mm,inner sep = 0pt ]
		\node (n0) at (0,0)  {\includegraphics[scale=\sclbm]{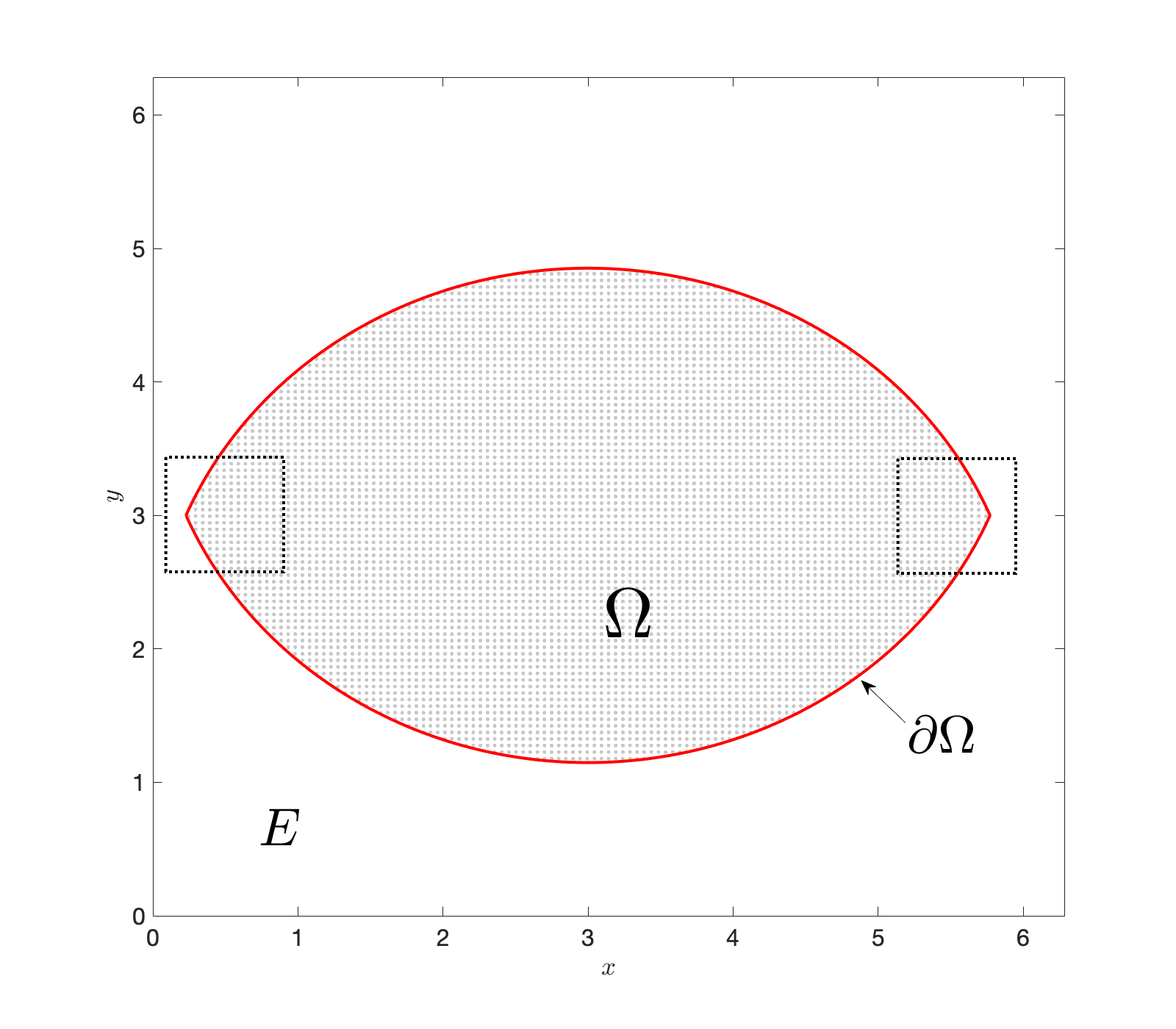}};
		\node (n2) at (-1.9,2.6)
		{\includegraphics[scale=\sclim]{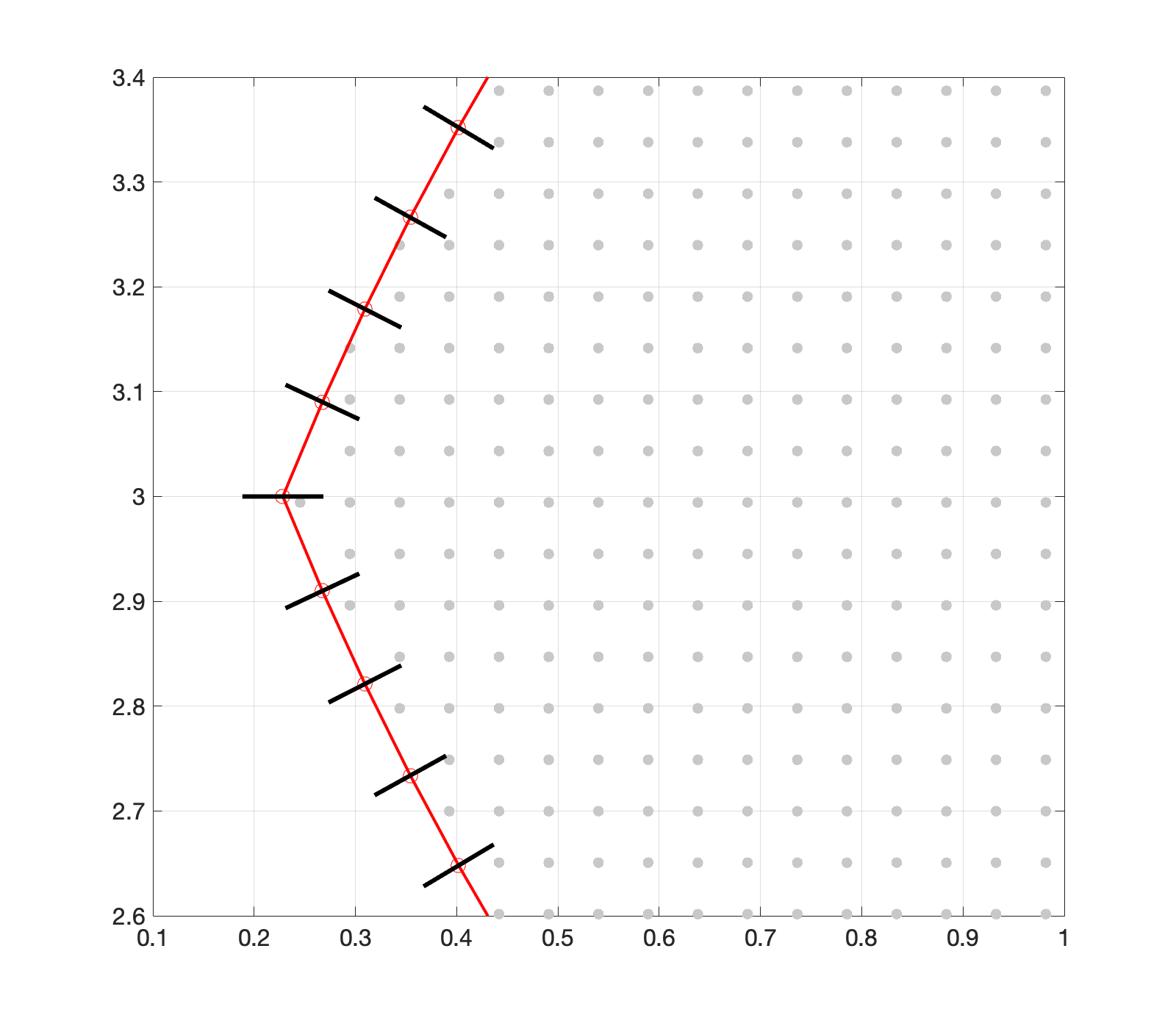}};
		\node (n3) at (2.8,2.6)
		{\includegraphics[scale=\sclim]{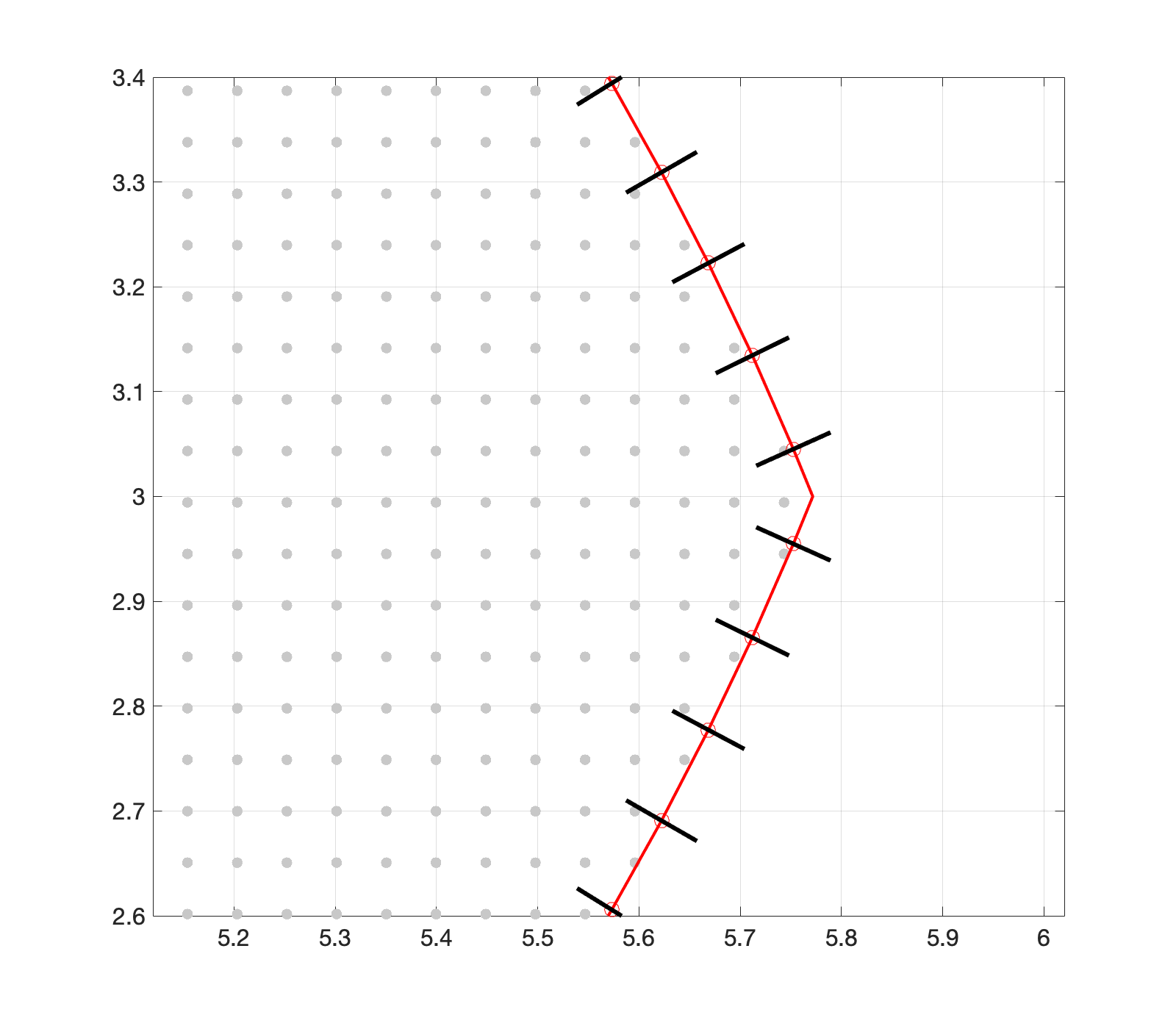}};        
		\node (n4) at (-1.5,3.1)  {(I)};
		\node (n5) at (2.4,3.1)  {(II)};
		\draw[densely dashed]   (-2.9,0.45) -- (-2.7,3.4) 
		(-2.1,-0.4) -- (-1.0,1.88);
		\draw[densely dashed]   (2.95,-0.4) -- (3.68,1.88) 
		(2.13,-0.4) -- (2.0,1.88);                 
		\end{tikzpicture}
		\label{eyeshape}
	}
	\subfigure[]
	{\begin{tikzpicture}
		[node distance = 10mm,inner sep = 0pt ]
		\node (n0) at (0,0)  {\includegraphics[scale=\sclbm]{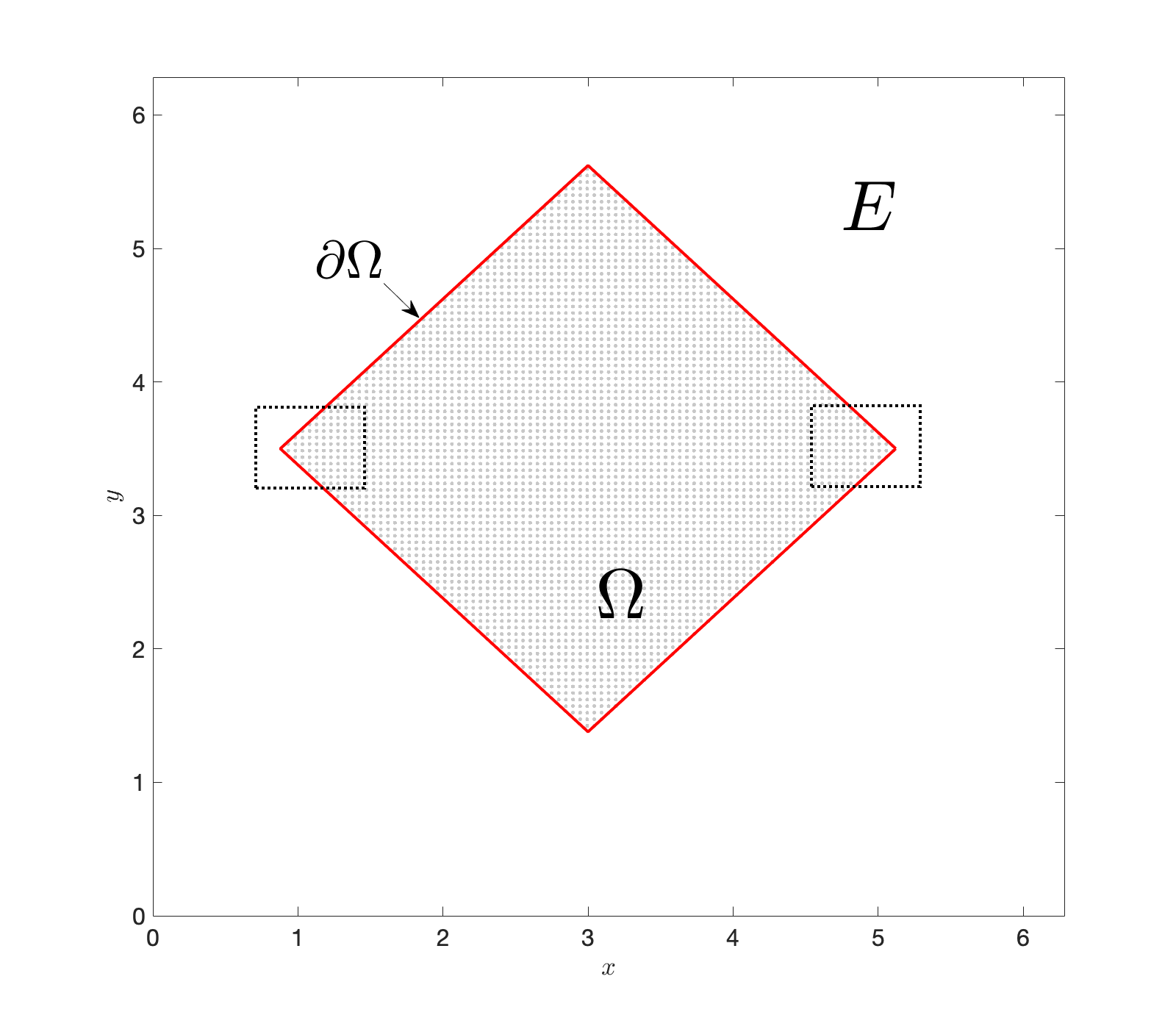}};
		\node (n2) at (-1.84,-1.75)
		{\includegraphics[scale=\sclim]{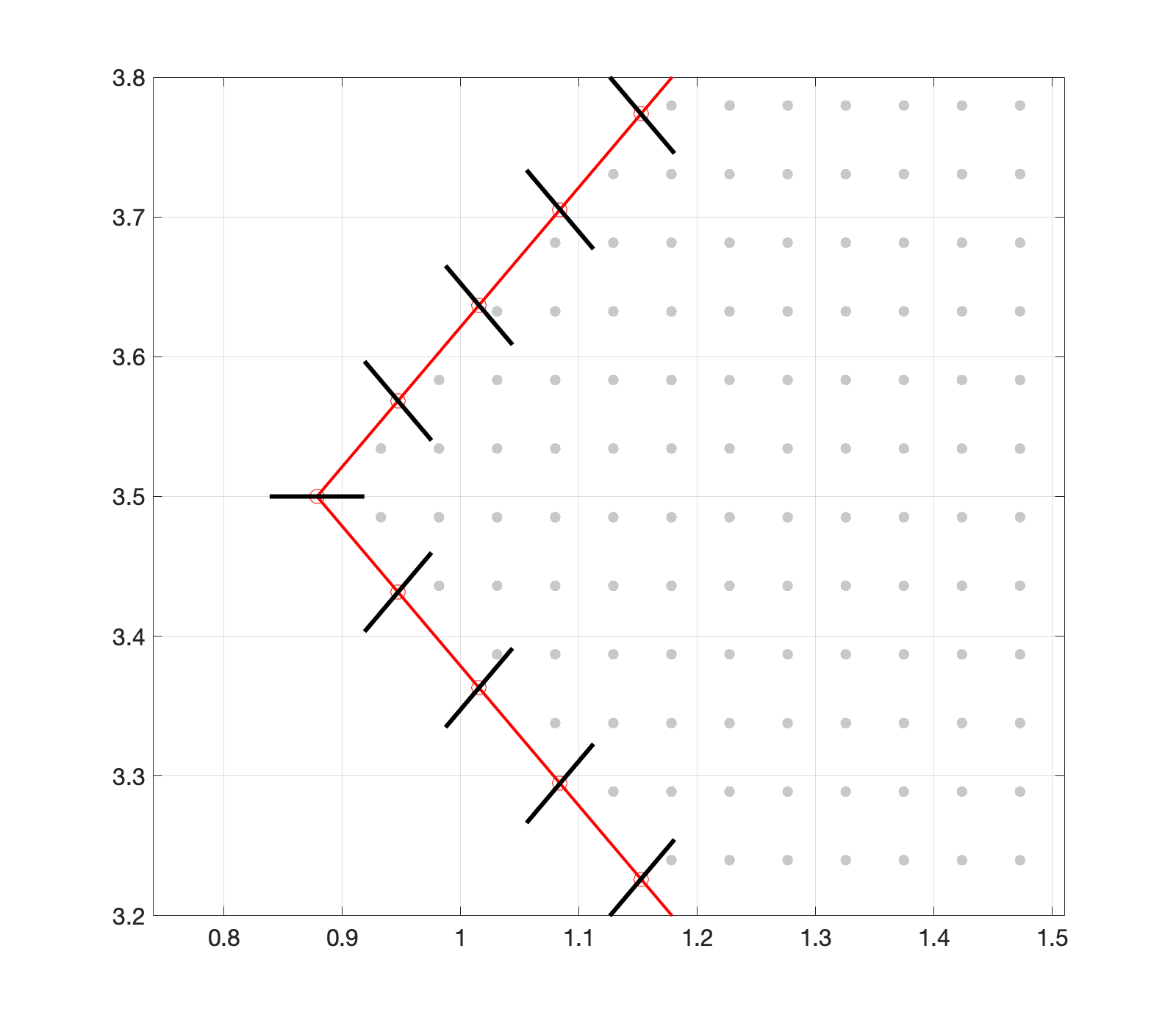}};
		\node (n3) at (2,-1.75)
		{\includegraphics[scale=\sclim]{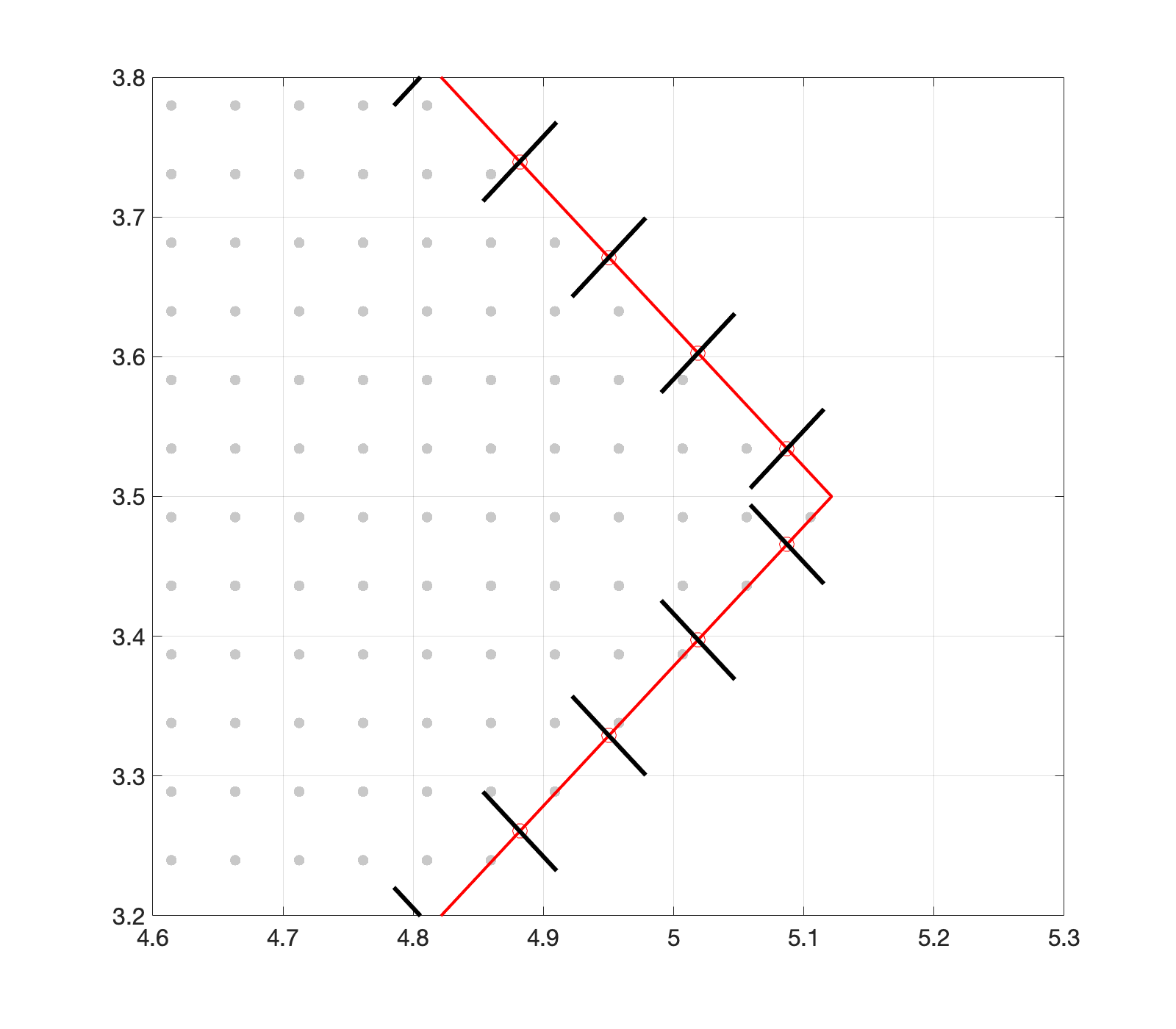}};     
		\node (n4) at (-1.25,-2.25)  {(I)};
		\node (n5) at (1.5,-2.25)  {(II)};   
		\draw[densely dashed]   (-2.28,0.72) -- (-2.62,-0.9) 
		(-1.52,0.72) -- (-1.0,-0.9);
		\draw[densely dashed]   (2.3,0.72) -- (2.87,-0.9) 
		(1.54,0.72) -- (1.2,-0.9);                 
		\end{tikzpicture}
		\label{diamshape}
	}
	\caption{(a) The eye-shaped domain $\Omega$ along with its boundary $\partial \Omega$ and extension region $E$. The fixed grid is shown in $\Omega$, along with the boundary nodes, for $N = 2^7$ and $n_\text{b} = 144$. In discretization approach (I), we place nodes on the corners and choose normal vectors symmetrically. In (II), we avoid placing nodes on the corners altogether. (b) The diamond domain $\Omega$ with the fixed grid and boundary nodes for $N = 2^7$ and $n_\text{b} = 124$. Insets (I) and (II) illustrate the two boundary discretization approaches.}\label{cornerdoms}
\end{figure}
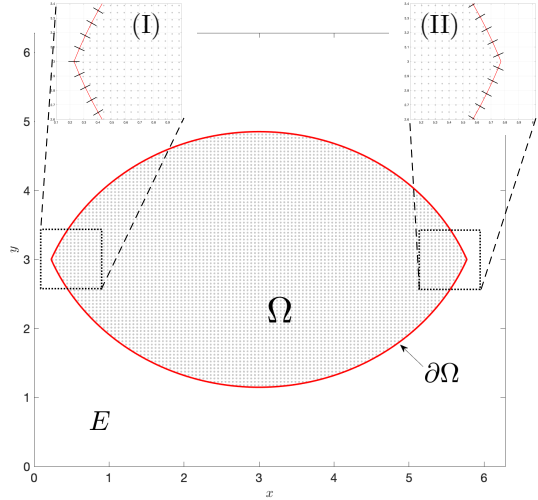
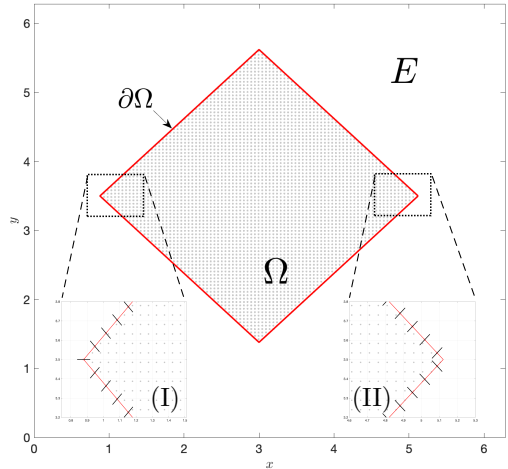

\begin{figure}[tbph]
	\centering
	\subfigure[]
	{\includegraphics[scale=\sclbm]{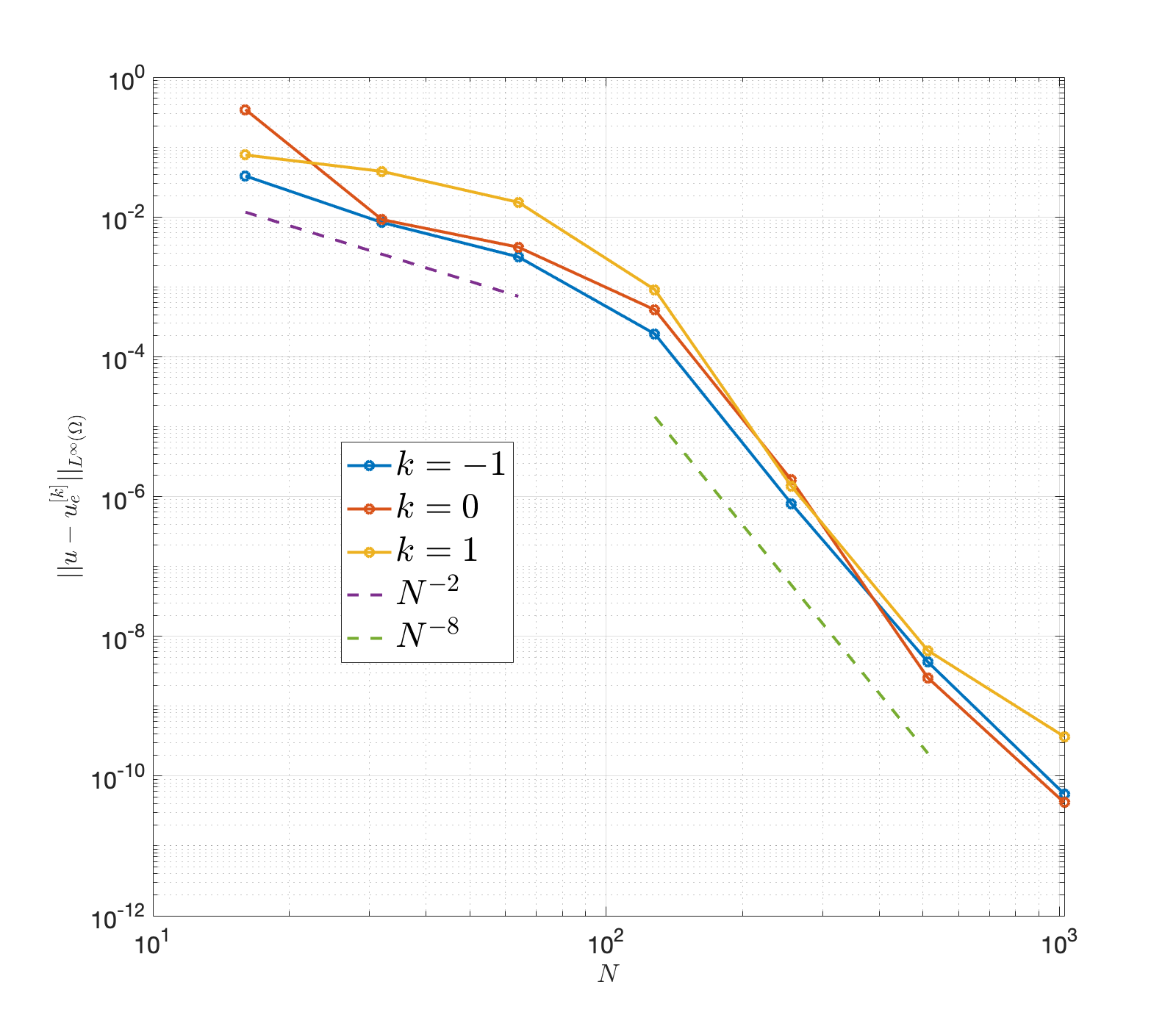}
		\label{lap2deye1}
	}
	\subfigure[]
	{\includegraphics[scale=\sclbm]{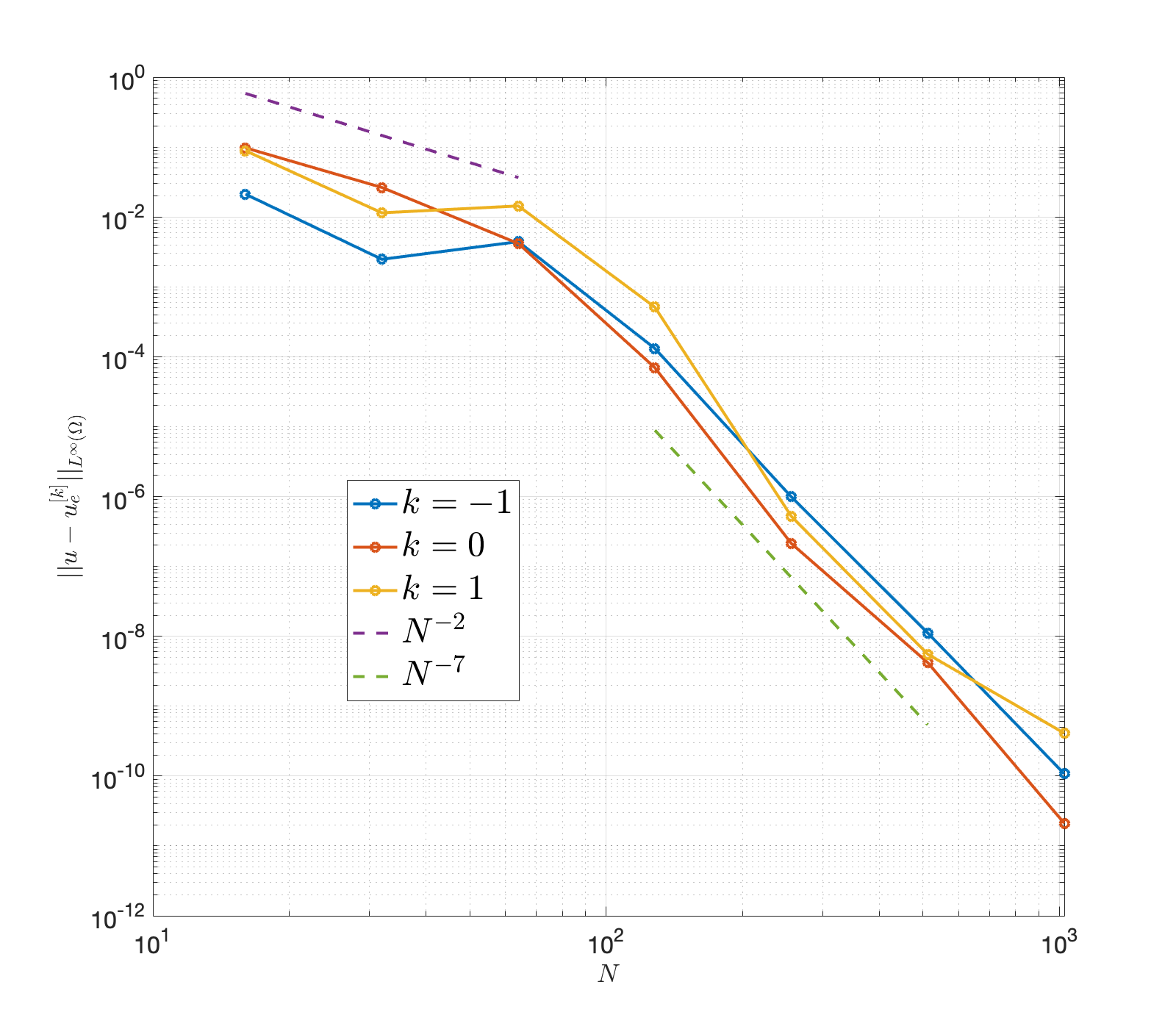}
		\label{lap2deye0}
	}
	\caption{(a) The convergence results for the Poisson equation solved on the eye-shaped domain shown in Figure \ref{eyeshape} with boundary nodes chosen as in (I). The solutions appear to converge sub-geometrically again, as shown by the speed-up in error decay rates from $O(N^{-2})$ to $O(N^{-8})$ for increasing $N$. Moreover, there appears to be little improvement for increasing values of $k$. (b) The same problem and domain as in (a) with boundary nodes as shown in inset (II) of Figure \ref{eyeshape}. The convergence properties are broadly similar to those in (a). However, the errors somewhat fluctuate for small $N$ before decaying at steady rates that are marginally slower than in (a).}\label{lap2deye}
\end{figure}

\begin{figure}[tbph]
	\centering
	\subfigure[]
	{\includegraphics[scale=\sclbm]{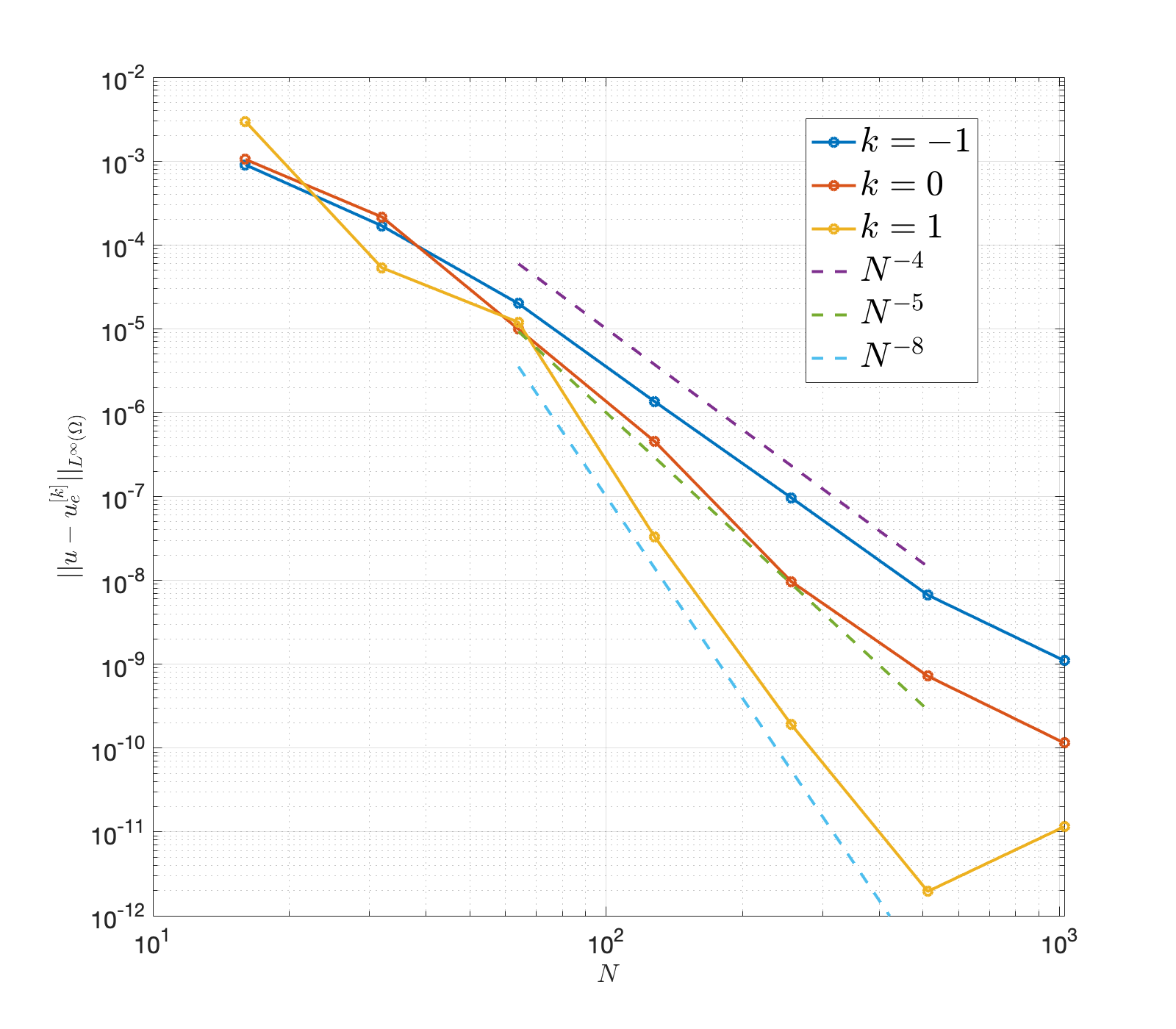}
		\label{lap2ddiam1}
	}
	\subfigure[]
	{\includegraphics[scale=\sclbm]{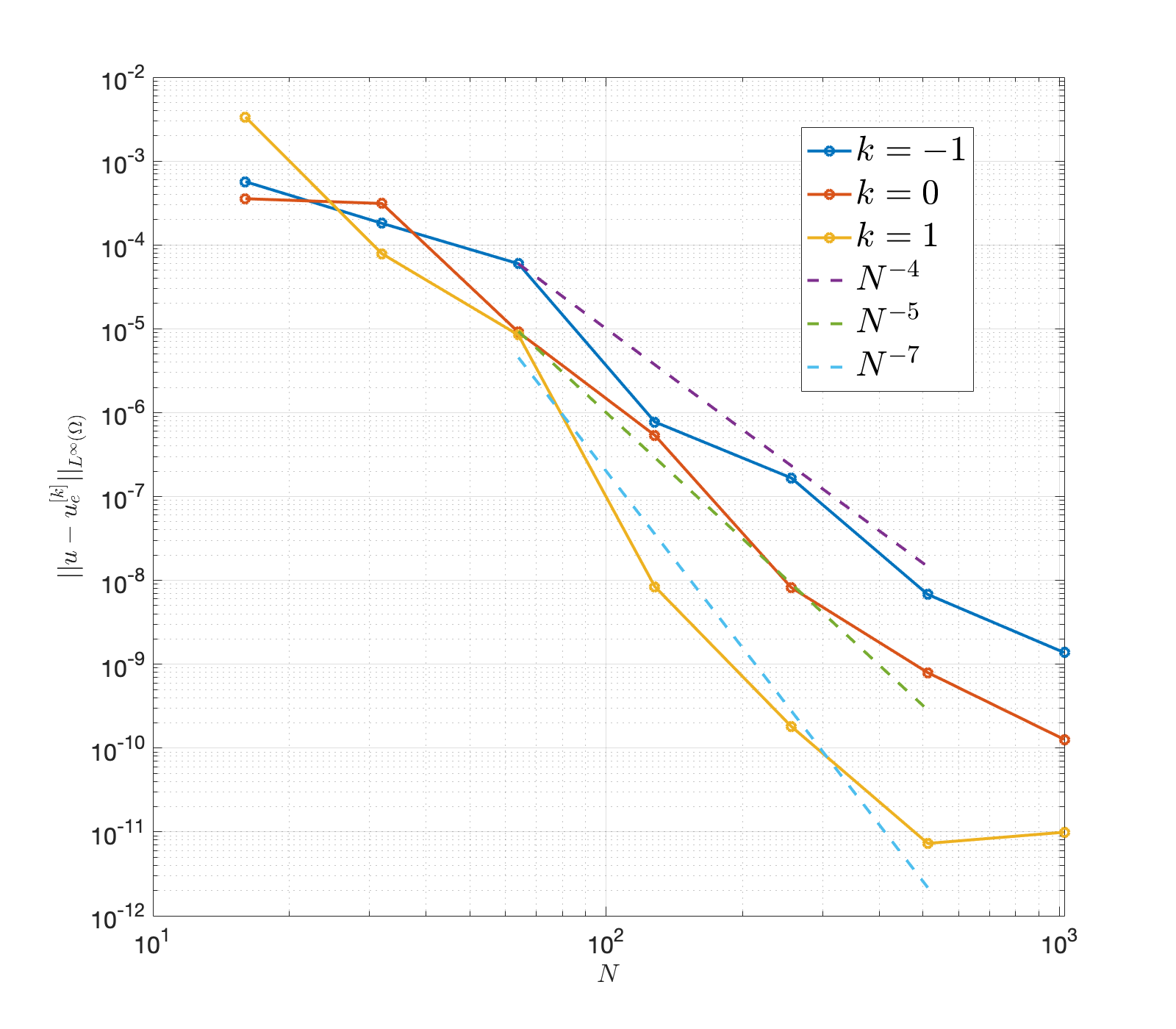}
		\label{lap2ddiam0}
	}
	\caption{(a) The convergence results for the Poisson equation solved on the diamond domain shown in Figure \ref{diamshape} with boundary nodes and normal vectors chosen as in (I). Unlike the earlier Poisson tests in 2D (Figures \ref{lap2d0} and \ref{lap2deye}), the errors converge algebraically in $1/N$. The convergence rates are still much higher than the expected $O(N^{-(k+3)})$. Moreover, we observe significant improvement as the value of $k$ is increased. (b) The same problem and domain as in (a), with the nodes placed away from the corners, as in inset (II) of Figure \ref{diamshape}. While the errors decay are similar to those seen in (a), they fluctuate more. For $k = 1$, the asymptotic rate is also slightly slower.}\label{lap2ddiam}
\end{figure}

To investigate this further, we next solve the Poisson equation on more challenging domains. Traditional fixed grid methods work reasonably well for domains with smooth boundaries, as in (\ref{lap2dprob1}), but suffer from poor performance when applied to non-smooth boundaries.

First, we consider the eye-shaped domain shown in Figure \ref{eyeshape}. This is centered at $(3,3)$ and built out of a pair of arcs, each subtended by an angle of $\Theta = 3\pi/4$ and radius $R = 3$. As it is not immediately obvious how the corners should be treated while discretizing the boundary, we explore two approaches. First, as shown in inset (I), we place nodes on the corners and choose a symmetrical normal vector direction. For the second approach, we avoid placing nodes on the corners altogether and jump over from one arc to the other, as shown in inset (II). In both instances, we use $n_\text{b} = \left\lceil \frac{R\Theta N}{2\pi} \right\rceil$ to ensure, as earlier, that boundary node spacing is roughly twice the grid spacing. We also solve the problem over the diamond domain shown in Figure \ref{diamshape}. The boundary in this case is a square of side-length $s = 3$ centered at $(3,3.5)$ and rotated by $45^\circ$. We take $n_\text{b} = \left\lceil \frac{s N}{\pi} \right\rceil$ and again consider both boundary discretization approaches (insets (I) and (II)). The number of extension functions are calculated by (\ref{Jchoice})

We test our algorithm on these domains by solving a problem with a known solution. We take $u(x,y) = \frac{1}{x^2 + y^2}$ with corresponding forcing $f(x,y) = -\frac{4}{(x^2+y^2)^2}$. Note that both the solution and the forcing possess singularities outside the physical domains, making this another instance of the ``mountain-in-fog'' test. We enforce regularity constraints of type (\ref{Tkeq2}) to get around the singularities and non-periodicity of the natural extension of $f$ on $\mathbb{T}^2$.

The errors for the eye-shaped domain, shown in Figure \ref{lap2deye}, can be seen to decay sub-geometrically again. The convergence for node placement of type (I) is better than for type (II), which fluctuates more and is asymptotically slower. In either case, the convergence rates comfortably exceed the expected $O(N^{-(k+3)})$. However, increasing the value of $k$ does not appear to lead to improved results and may even lead to marginally larger errors.

This is not the case for the diamond domain. The error plots, shown in Figure \ref{lap2ddiam}, show that the convergence improves significantly with increasing $k$. In this case, the errors can be seen to decay algebraically in $1/N$, although still exceeding the $O(N^{-(k+3)})$ threshold. As earlier, the performance is less prone to fluctuations when nodes are placed on corners. To sum up, these tests show that our technique works at least as well as designed for $C^0$ boundaries with more consistent performance with nodes placed on corners.

\subsection{Computing Eigenvalues on Arbitrary Domains}
Next, we use our method to find the eigenvalues of the Laplacian on the domains considered earlier. Since our technique allows for rapid and accurate inversion of elliptic operators, it is well-suited to the power method. This shall also serve as a useful test of the stability of our algorithm, in that it avoids spurious eigenvalues, which is a critical issue while solving time-dependent problems. Consider the eigenvalue problem,
\eqn{
	\left\{\begin{matrix}
		-\Delta u = \lambda u, & \text{on } \Omega, \\
		u = 0, & \text{on } \partial \Omega. 
	\end{matrix}\right.	\label{EigProb}
}
It is well-known that the eigenvalues of the Laplacian are real, positive and can be arranged as
\eqn{
	0 < \lambda_1 < \lambda_2 < \lambda_3 < \hdots	\label{eigseq}
}
To compute $\lambda_i$ we choose a real shift $\sigma$ such that
\eqn{
	|\lambda_i - \sigma| < |\lambda_j - \sigma|, \qquad j \neq i, \label{sigcond}	
}
so that the smallest eigenvalue (in the absolute sense) of $(-\Delta - \sigma)$ is $(\lambda_i - \sigma)$. Thus, applying the power method to $(-\Delta - \sigma)^{-1}$ should allow us to find $\lambda_i$. More precisely, we compute 
\eqn{
	\left\{\begin{matrix}
		(-\Delta - \sigma)v^{n+1} = u^n, & \text{on } \Omega, \\
		v^{n+1} = 0, & \text{on } \partial \Omega, \\
		u^{n+1} = v^{n+1}/\norm{v^{n+1}}_{L^2(\Omega)}, & \\
	\end{matrix}\right.	\label{EigProbMod}
}
and $\tilde{\lambda}_{n+1} = \ip{u^{n+1},-\Delta u^{n+1}}$, for $n \geq 0$. The initial seed $u^0$ is chosen randomly; the values $\{\tilde\lambda_{n}\}$ then converge geometrically to the desired eigenvalue $\lambda_i$ with high probability. 

Our approach to solving (\ref{EigProbMod}) requires $\Omega$ to be embedded in the computational domain $C = \mathbb{T}^d$. To avoid imposing an averaging condition of the form (\ref{Avgeq}), we choose $\sigma$ so that $(-\Delta - \sigma)$ is invertible on $C$; this can be achieved easily by choosing $\sigma$ to be non-integer since $(-\Delta - \sigma)$ fails to be invertible on $C$ if and only if $\sigma = \sum_{l = 1}^d m_l^2$, with $m_l \in \mathbb{Z}$, for all $l$. By varying $\sigma$, we can find all the eigenvalues.

\begin{table}	
	\centering
	\begin{tabular}{ccc} 
		\hline\hline
		\textbf{Figure \ref{discshape}} & \textbf{Figure \ref{eyeshape}} & \textbf{Figure \ref{diamshape}}  \\ 
		\hline
		$\lambda_i$ & $\lambda_i|\Omega|$ & $\lambda_i|\Omega|$ \\
		\hline
		$0.219308$ & $19.3222$ & $19.7392$ \\
		$1.06247$ & $41.0926$ & $49.3480$ \\
		$1.23047$ & $56.7845$ &  $78.9568$ \\
		$1.63043$ & $71.3660$ & $98.6960$ \\
		$2.19414$ & $89.1973$ & $128.305$ \\
		$2.55057$ & $110.091$ & $167.783$ \\
		$3.75893$ & $116.318$ & $177.653$ \\
		\hline
	\end{tabular}\caption{The first seven eigenvalues of $(-\Delta)$ for the various two-dimensional domains considered earlier, with homogeneous Dirichlet boundary conditions, shown up to six significant figures. For the eye-shape and diamond, they have been scaled by the areas to make them independent of the parameters used to define these domains.} \label{evaltab}
\end{table}

The iterations (\ref{EigProbMod}) are continued until the deviations 
\eqn{
	d_n = \max\{|\tilde {\lambda}_{n+1} - \tilde {\lambda}_n|, \norm{u^{n+1} - u^n}_{L^2(\Omega)}\} \label{devs}
} 
fall below a pre-determined tolerance $\tau$. In our computations, we use $\tau = 10^{-10}$ with $k = -1$ for the extension and $N = 2^9$ for the grid-size. 

The first few computed eigenvalues, over the domains considered earlier, are shown in Table \ref{evaltab}. The eigenvalues for the ``interior'' problems in Figures \ref{eyeshape} and \ref{diamshape} have been scaled by the areas of the domains, to make the results independent of the side-lengths or radii used. We use boundary discretizations of type (I) for these domains because, as established in Figures \ref{lap2deye} and \ref{lap2ddiam}, this choice leads to more accurate results. The accuracy of these computations can be assessed by comparing them to analytically calculated values. For instance, the eigenvalues for the diamond in Figure \ref{diamshape} obey
\eqn{
	\lambda_{m,n} = \frac{\pi^2}{|\Omega|}(m^2 + n^2), \qquad m,n \in \mathbb{Z}_+. \label{diameigs}	
} 
It can be seen that the corresponding values in Table \ref{evaltab} are indeed just the appropriate multiples of $\pi^2$.

\subsection{Heat Equation in Two Dimensions}
Finally, we apply the methodology to the heat equation. As discussed earlier, applying the BDF-4 time discretization leads to a Helmholtz problem of type (\ref{bdf41}) at each time step. We consider the ``external'' problem on the domain in Figure \ref{discshape} and calculate the forcing, initial condition and boundary values from the exact solution
\eqn{
	u(t,x,y) = e^{\sin(x)}\cos(y)\cos(t). \label{heat2Dex}	
}  

\begin{figure}[tbph]
	\centering
	{\includegraphics[scale = \sclas]{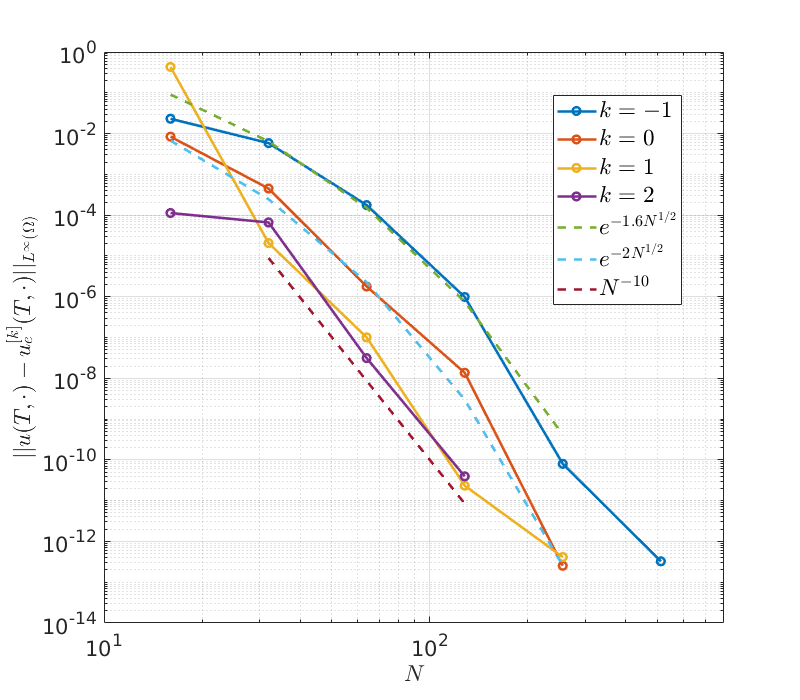}}
	\caption{The $L^{\infty}$ errors at final time $T = 2$ in the solutions to the 2D heat equation. The time-step sizes are chosen sufficiently small so that time integration errors are negligible. The convergence in space is sub-geometric for $k \leq 0$, as seen earlier on a domain with a smooth boundary. For smoother extensions, the convergence still comfortably exceeds expectations.} \label{heat2dinf}
	
\end{figure}

As we are using an exact solution, we do not need to jump-start the multi-step scheme by using a single step method for the initial few steps. As in the time-dependent example seen earlier, we opt for regularity constraints of the form (\ref{Tkeq2}). In contrast, however, we vary the time-step size with grid spacing due to stability considerations. Recall that we invert the Helmholtz operator $\hat{\mathcal{L}} =  \mathbb{I} - \frac{12\Delta t}{25} \Delta $ at each time-step; an exceedingly small time-step makes it harder to damp out the high frequencies. Instead, we find it more beneficial to scale the number of time steps with the number of grid points. Specifically, we opt for $\Delta t = 1/4N$, so that $\Delta t = \Delta x/8\pi$. The high order time integration scheme then ensures that the time-stepping errors are negligible, allowing us to assess the accuracy in space of the solutions to the iterated Helmholtz problems.

The solutions to an elliptic problem on a domain with a smooth boundary have already been seen to converge faster than $O(N^{-(k+3)})$ and, in particular, sub-geometrically for $k = -1$. In addition, that repeated iterations of the solver lead to accurate results has been established by the solutions to the eigenvalue problems. These ingredients combine to yield the $L^{\infty}$ convergence plots shown in Figure \ref{heat2dinf}. We note that both $k = -1$ and $k = 0$ extensions appear to yield sub-geometric convergence. For higher $k$, the errors still decay faster than expected.

\section{Conclusion}\label{conc}

In this study, we have introduced a technique for solving elliptic problems on arbitrary domains. Our approach uses and further develops the ideas and insights that power fixed Cartesian grid techniques such as the IB and the IBSE methods. In this sense, it may be seen as a next step in the sequence. At the same time, it eschews many of the tools that are pervasive in these approaches, such as discretized delta functions and local corrections via spreading operators, in the process making it more accurate. 

A signature feature of our methodology is that one can obtain arbitrary orders of accuracy by appropriately setting the regularity of the extension to the forcing. The manner in which this is done avoids the ``mountain-in-fog'' problem and enables one to solve problems whose analytic solutions may have ill-behaved natural extensions. The use of NUFFT algorithms for interpolation, apart from speeding up the computations, also ensures that there is no barrier to the highest achievable accuracy. Moreover, the technique is demonstrably stable: repeatedly iterating the solvers, as we did while solving the time-dependent problems and computing the eigenvalues, does not lead to numerical blow-up, to which some spectral methods are susceptible \cite{boyd2001chebyshev,dawkins1998origin,gottlieb1977numerical}.

For one-dimensional problems, the observed rates of convergence are in perfect agreement with theory, indicating the soundness of our approach. The performance of our technique for two-dimensional domains, however, is much better than anticipated. The convergence rates comfortably exceed the expected rates and, in specific cases, appear to be sub-geometric.

The simplicity of this approach makes it easy to extend it to higher dimensions. Another avenue for exploration is in the development of algorithms for solving fluid equations, such as the Stokes and Navier--Stokes models. A further potentially fruitful extension is to models of viscoelastic fluids. Low order methods generally fail to capture the stress values close to boundaries, which limits their usefulness in such regimes \cite{stein2019convergent}. Since our technique allows the order of accuracy to be set arbitrarily, its application to these problems has the potential to lead to significant advances.       

\section{Acknowledgements}
This work is supported by the National Science Foundation awards DMS 1664645, OAC 1450327, OAC 1652541 and OAC 1931516. The authors thank Robert D. Guy, Becca Thomases, M. Gregory Forest, J. Thomas Beale, Aleksandar Donev, Ebrahim M. Kolahdouz, and Charles Puelz for their insightful suggestions and comments. 

\bibliography{refs}

\begin{thebibliography}{10}
\expandafter\ifx\csname url\endcsname\relax
  \def\url#1{\texttt{#1}}\fi
\expandafter\ifx\csname urlprefix\endcsname\relax\def\urlprefix{URL }\fi
\expandafter\ifx\csname href\endcsname\relax
  \def\href#1#2{#2} \def\path#1{#1}\fi

\bibitem{shewchuk2002good}
J.~R. Shewchuk, What is a good linear finite element? interpolation,
  conditioning, anisotropy, and quality measures (preprint), University of
  California at Berkeley 73 (2002) 137.

\bibitem{helsing2008evaluation}
J.~Helsing, R.~Ojala, On the evaluation of layer potentials close to their
  sources, Journal of Computational Physics 227~(5) (2008) 2899--2921.

\bibitem{beale2001method}
J.~T. Beale, M.-C. Lai, A method for computing nearly singular integrals, SIAM
  Journal on Numerical Analysis 38~(6) (2001) 1902--1925.

\bibitem{fryklund2019integral}
F.~Fryklund, M.~C.~A. Kropinski, A.-K. Tornberg, An integral equation based
  numerical method for the forced heat equation on complex domains, arXiv
  preprint arXiv:1907.08537.

\bibitem{peskin1972flow}
C.~S. Peskin, Flow patterns around heart valves: a numerical method, Journal of
  computational physics 10~(2) (1972) 252--271.

\bibitem{peskin1977numerical}
C.~S. Peskin, Numerical analysis of blood flow in the heart, Journal of
  computational physics 25~(3) (1977) 220--252.

\bibitem{griffith2007adaptive}
B.~E. Griffith, R.~D. Hornung, D.~M. McQueen, C.~S. Peskin, An adaptive,
  formally second order accurate version of the immersed boundary method,
  Journal of computational physics 223~(1) (2007) 10--49.

\bibitem{bhalla2013unified}
A.~P.~S. Bhalla, R.~Bale, B.~E. Griffith, N.~A. Patankar, A unified
  mathematical framework and an adaptive numerical method for fluid--structure
  interaction with rigid, deforming, and elastic bodies, Journal of
  Computational Physics 250 (2013) 446--476.

\bibitem{griffith2012immersed}
B.~E. Griffith, Immersed boundary model of aortic heart valve dynamics with
  physiological driving and loading conditions, International Journal for
  Numerical Methods in Biomedical Engineering 28~(3) (2012) 317--345.

\bibitem{huang2012three}
W.-X. Huang, C.~B. Chang, H.~J. Sung, Three-dimensional simulation of elastic
  capsules in shear flow by the penalty immersed boundary method, Journal of
  Computational Physics 231~(8) (2012) 3340--3364.

\bibitem{kou2015fully}
W.~Kou, A.~P.~S. Bhalla, B.~E. Griffith, J.~E. Pandolfino, P.~J. Kahrilas,
  N.~A. Patankar, A fully resolved active musculo-mechanical model for
  esophageal transport, Journal of computational physics 298 (2015) 446--465.

\bibitem{seol2016immersed}
Y.~Seol, W.-F. Hu, Y.~Kim, M.-C. Lai, An immersed boundary method for
  simulating vesicle dynamics in three dimensions, Journal of Computational
  Physics 322 (2016) 125--141.

\bibitem{stein2016immersed}
D.~B. Stein, R.~D. Guy, B.~Thomases, {Immersed boundary smooth extension: a
  high-order method for solving PDE on arbitrary smooth domains using Fourier
  spectral methods}, Journal of Computational Physics 304 (2016) 252--274.

\bibitem{stein2017immersed}
D.~B. Stein, R.~D. Guy, B.~Thomases, {Immersed Boundary Smooth Extension
  (IBSE): A high-order method for solving incompressible flows in arbitrary
  smooth domains}, Journal of Computational Physics 335 (2017) 155--178.

\bibitem{boyd2002comparison}
J.~P. Boyd, {A comparison of numerical algorithms for Fourier extension of the
  first, second, and third kinds}, Journal of Computational Physics 178~(1)
  (2002) 118--160.

\bibitem{trefethen2013approximation}
L.~N. Trefethen, Approximation theory and approximation practice, Vol. 128,
  SIAM, 2013.

\bibitem{clenshaw1955note}
C.~W. Clenshaw, {A note on the summation of Chebyshev series}, Mathematics of
  Computation 9~(51) (1955) 118--120.

\bibitem{leveque1994immersed}
R.~J. Leveque, Z.~Li, {The immersed interface method for elliptic equations
  with discontinuous coefficients and singular sources}, SIAM Journal on
  Numerical Analysis 31~(4) (1994) 1019--1044.

\bibitem{li2006immersed}
Z.~Li, K.~Ito, The immersed interface method: numerical solutions of PDEs
  involving interfaces and irregular domains, SIAM, 2006.

\bibitem{fedkiw1999non}
R.~P. Fedkiw, T.~Aslam, B.~Merriman, S.~Osher, et~al., {A non-oscillatory
  Eulerian approach to interfaces in multimaterial flows (the ghost fluid
  method)}, Journal of computational physics 152~(2) (1999) 457--492.

\bibitem{marques2011correction}
A.~N. Marques, J.-C. Nave, R.~R. Rosales, {A correction function method for
  Poisson problems with interface jump conditions}, Journal of Computational
  Physics 230~(20) (2011) 7567--7597.

\bibitem{shirokoff2015sharp}
D.~Shirokoff, J.-C. Nave, {A sharp-interface active penalty method for the
  incompressible Navier--Stokes equations}, Journal of Scientific Computing
  62~(1) (2015) 53--77.

\bibitem{bruno2010high}
O.~P. Bruno, M.~Lyon, {High-order unconditionally stable FC-AD solvers for
  general smooth domains I. Basic elements}, Journal of Computational Physics
  229~(6) (2010) 2009--2033.

\bibitem{lyon2010high}
M.~Lyon, O.~P. Bruno, {High-order unconditionally stable FC-AD solvers for
  general smooth domains II. Elliptic, parabolic and hyperbolic PDEs;
  theoretical considerations}, Journal of Computational Physics 229~(9) (2010)
  3358--3381.

\bibitem{agress2018novel}
D.~Agress, P.~Guidotti, A novel optimization approach to fictitious domain
  methods, arXiv preprint arXiv:1808.02158.

\bibitem{agress2019smooth}
D.~Agress, P.~Guidotti, D.~Yan, {The Smooth Selection Embedding Method with
  Chebyshev Polynomials}, arXiv preprint arXiv:1902.03713.

\bibitem{fryklund2018partition}
F.~Fryklund, E.~Lehto, A.-K. Tornberg, Partition of unity extension of
  functions on complex domains, Journal of Computational Physics 375 (2018)
  57--79.

\bibitem{klinteberg2019fast}
L.~af~Klinteberg, T.~Askham, M.~C. Kropinski, {A fast integral equation method
  for the two-dimensional Navier-Stokes equations}, Journal of Computational
  Physics 409 (2020) 109353.

\bibitem{evans2010partial}
L.~C. Evans, Partial differential equations, Vol.~19, American Mathematical
  Soc., 2010.

\bibitem{greengard2004accelerating}
L.~Greengard, J.-Y. Lee, {Accelerating the nonuniform fast Fourier transform},
  SIAM review 46~(3) (2004) 443--454.

\bibitem{lee2005type}
J.-Y. Lee, L.~Greengard, {The type 3 nonuniform FFT and its applications},
  Journal of Computational Physics 206~(1) (2005) 1--5.

\bibitem{demmel1997applied}
J.~W. Demmel, Applied numerical linear algebra, SIAM, 1997.

\bibitem{kallemov2016immersed}
B.~Kallemov, A.~Bhalla, B.~Griffith, A.~Donev, An immersed boundary method for
  rigid bodies, Communications in Applied Mathematics and Computational Science
  11~(1) (2016) 79--141.

\bibitem{boyd2001chebyshev}
J.~P. Boyd, Chebyshev and Fourier spectral methods, Courier Corporation, 2001.

\bibitem{dawkins1998origin}
P.~T. Dawkins, S.~R. Dunbar, R.~W. Douglass, The origin and nature of spurious
  eigenvalues in the spectral tau method, Journal of Computational Physics
  147~(2) (1998) 441--462.

\bibitem{gottlieb1977numerical}
D.~Gottlieb, S.~A. Orszag, Numerical analysis of spectral methods: theory and
  applications, Vol.~26, Siam, 1977.

\bibitem{stein2019convergent}
D.~B. Stein, R.~D. Guy, B.~Thomases, {Convergent solutions of Stokes Oldroyd-B
  boundary value problems using the Immersed Boundary Smooth Extension (IBSE)
  method}, Journal of Non-Newtonian Fluid Mechanics 268 (2019) 56--65.

\end{thebibliography}

\end{document}